\documentclass[a4paper,leqno]{amsart}
\usepackage{amsfonts}
\theoremstyle{plain}
\begingroup

\newtheorem{theorem}{Theorem}[section]
\newtheorem{lemma}[theorem]{Lemma}
\newtheorem{proposition}[theorem]{Proposition}
\newtheorem{corollary}[theorem]{Corollary}
\newtheorem{remark}[theorem]{Remark}
\newtheorem{definition}[theorem]{Definition}
\endgroup

\theoremstyle{definition}
\theoremstyle{remark}

\mathchardef\emptyset="001F

\numberwithin{equation}{section}

\newcommand{\e}{\varepsilon}
\newcommand{\Om}{\Omega}
\newcommand{\Omk}{\Om\setmeno K}

\newcommand{\R}{{\mathbb R}}
\newcommand{\wto}{{\rightharpoonup}}
\newcommand{\setmeno}{\!\setminus\!}
\newcommand{\pinfty}{{+}\infty}

\newcommand{\huno}{{\mathcal H}^{1}}
\newcommand{\E}{{\mathcal E}}
\newcommand{\hb}{L^{1,2}}
\newcommand{\gdot}{{\dot{g}}}
\newcommand{\K}{{\mathcal K}(\overline\Om)}
\newcommand{\Kf}{{\mathcal K}^f(\overline\Om)}
\newcommand{\Kp}{{\mathcal K}_p(\overline\Om)}
\newcommand{\Kpf}{{\mathcal K}_p^f(\overline\Om)}
\newcommand{\Km}{{\mathcal K}_m(\overline\Om)}
\newcommand{\Kmf}{{\mathcal K}_m^f(\overline\Om)}
\newcommand{\Kml}{{\mathcal K}_m^\lambda(\overline\Om)}
\newcommand{\Kone}{{\mathcal K}_1(\overline\Om)}
\newcommand{\Konef}{{\mathcal K}_1^f(\overline\Om)}
\newcommand{\Konel}{{\mathcal K}_1^\lambda(\overline\Om)}
\newcommand{\KA}{{\mathcal K}_{q+m-h}^f(\overline A)}

\title[Quasi-static growth of brittle fractures]
{A  model for the quasi-static growth\\
of brittle fractures:\\ existence
and approximation results
  }
\author[Gianni Dal Maso]{Gianni Dal Maso}
\address[Gianni Dal Maso]{SISSA, Via Beirut 2-4, 34014 Trieste,
Italy}
\email[Gianni Dal Maso]{dalmaso@sissa.it}
\author[Rodica Toader]{Rodica Toader}
\address[Rodica Toader]{SISSA, Via Beirut 2-4, 34014 Trieste,
Italy}
\email[Rodica Toader]{toader@sissa.it}
\begin{document}

\begin{abstract}
We give a precise mathematical formulation of a variational model for
the irreversible quasi-static evolution of brittle fractures proposed by
G.A.~Francfort and J.-J.~Marigo, and based on Griffith's theory of
crack growth.
In the two-dimensional case we prove an existence result for the
quasi-static evolution and show that the total energy is an absolutely
continuous function of time, although we can not exclude that
the bulk energy and the surface
energy may present some jump discontinuities. This existence result
is proved by a time discretization process, where at each step a
global energy minimization is performed, with the constraint that the
new crack contains all cracks formed at the previous time steps. This
procedure provides an effective way to approximate the
continuous time evolution.
\end{abstract}
\maketitle
{\small

\bigskip


\bigskip
}
\bigskip
\bigskip

\begin{section}{INTRODUCTION}

Since the pioneering work of A.~Griffith \cite{Gri}, the growth of a
brittle fracture is considered to be the result of the competition
between the energy spent to increase the crack and the corresponding
release of bulk energy. This idea is the basis of the celebrated
Griffith's criterion for crack growth (see, e.g., \cite{SL}), and is used to
study the crack propagation along a preassigned path.
The actual path followed by
the crack is often determined by using different criteria
(see, e.g., \cite{ES}, \cite{SL}, \cite{SM}).

Recently G.A.~Francfort and J.-J.~Marigo \cite{FraMar3} proposed
a variational model for the quasi-static growth of brittle
fractures, based on Griffith's theory,
where the interplay between bulk and surface energy
determines also the crack path.

The purpose of this paper is to give a precise mathematical
formulation of a variant of this model in the {\it two-dimensional
case\/},
and to prove an existence result for the {\it quasi-static evolution of a
fracture\/} by using the {\it time discretization method\/} proposed in
\cite{FraMar3}.

To simplify the mathematical description of the model, we consider
only {\it linearly elastic homogeneous isotropic
materials\/}, with Lam\'e coefficients
$\lambda$ and $\mu$. We restrict our analysis to
the case of an {\it anti-plane shear\/}, where
the reference configuration is an infinite cylinder ${\Om{\times}\R}$, with
$\Omega\subset \R^2$, and the displacement has the special form
$(0,0,u(x_1,x_2))$ for every $(x_1,x_2,y)\in {\Om{\times}\R}$.
We assume also that the cracks  have the form
$K{\times}\R$, where $K$ is a compact set in $\overline\Om$.
In this case the notions of bulk energy and surface energy
refer to a finite portion of the cylinder determined by
two cross sections separated by a unit distance.
The {\it bulk energy\/} is given by
\begin{equation}\label{bulk}
\frac{\mu}2\int_{\Om\setminus K} |\nabla u|^2dx\,,
\end{equation}
while the {\it surface energy\/} is given by
\begin{equation}\label{surface}
k\,\huno(K)\,,
\end{equation}
where $k$
is a constant which depends on the toughness of the material, and
$\huno$ is the {\it one-dimensional Hausdorff measure\/}, which
coincides with the ordinary length in case $K$ is a rectifiable arc. For
simplicity we take $\mu=2$ and $k=1$ in (\ref{bulk})
and (\ref{surface}).

We assume that $\Om$ is a {\it connected bounded open set\/}
with a {\it Lipschitz boundary\/} $\partial\Om$. As in
\cite{FraMar3}, we fix a subset $\partial_D\Om$ of $\partial\Om$, on
which we want to  prescribe a {\it Dirichlet boundary condition\/} for $u$. We
assume that $\partial_D\Om$ has a {\it finite number of connected
components\/}.

Given a function $g$ on $\partial_D\Om$,
we consider the boundary condition
$u=g$ on $\partial_D\Omk$. We can not prescribe a Dirichlet
boundary condition on $\partial_D\Om\cap K$, because
the boundary displacement is not
transmitted through the crack, if the crack touches the boundary. Assuming
that {\it the fracture is traction
free\/}
(and, in particular, without friction),
the displacement $u$ in $\Omk$ is obtained by
{\it minimizing (\ref{bulk}) under the boundary condition
$u=g$ on $\partial_D\Omk$\/}.
The {\it total energy\/} relative to the boundary displacement
$g$ and to the crack determined by $K$ is therefore
\begin{equation}\label{e0}
\E(g,K)=\min_v \Big\{\int_{\Om\setminus K}|\nabla
v|^2dx+\huno(K) : v=g \hbox{ on }\partial_D\Om\setmeno K \Big\}\,.
\end{equation}
As $K$ is not assumed to be smooth,
we have to be careful in the precise mathematical formulation
of this minimum problem, which is given
at the beginning of  Section~\ref{conjugate}. The corresponding
existence result is based on  some properties of the {\it Deny-Lions
spaces\/}, that are described in Section~\ref{notation}.

In the theory developed in \cite{FraMar3} a crack with finite
surface energy is any compact subset $K$ of $\overline\Om$ with
$\huno(K)<\pinfty$. For technical reasons, that will be explained
later, we propose a variant of this model, where we prescribe an a
priori bound on the number of connected components of the
cracks. Without this restriction, some
convergence arguments used in the proof of our existence result
are not justified
by the present development of
the mathematical theories related to this subject.

We now describe our model of {\it quasi-static irreversible
evolution of a fracture\/} under the action of a {\it time dependent
boundary displacement\/} $g(t)$, $0\le t\le 1$.
As usual, we assume that $g(t)$ can be extended to a function, still
denoted by $g(t)$, which belongs to the Sobolev space $H^1(\Omega)$.
In addition, we assume that the function $t\mapsto g(t)$ is {\it
absolutely continuous\/} from $[0,1]$ into $H^1(\Omega)$.
Given an integer $m\ge 1$, let $\Kmf$ be the set of all compact
subsets $K$ of $\overline\Om$ having at most $m$ connected components
and with $\huno(K)<\pinfty$.
Following
the ideas of \cite{FraMar3}, given an initial crack $K_0\in \Kmf$,
we look for an {\it increasing family\/}
$K(t)$, $0\le t\le 1$, of cracks in $\Kmf$,
such that for any time
$t\in(0,1]$ the crack $K(t)$ {\it minimizes the total energy\/}
$\E(g(t),K)$ among all cracks in $\Kmf$ which contain
all previous cracks $K(s)$, $s<t$. For $t=0$ we assume that
$K(0)$ minimizes $\E(g(0),K)$ among all  cracks in $\Kmf$
which contain $K_0$.

This minimality condition for every time $t$
is inspired by Griffith's analysis of the energy balance. The constraint given
by the presence of the previous
cracks reflects the {\it irreversibility of the evolution\/}
and the {\it absence of a healing process\/}. In addition to this minimality
condition we require also that
$\frac{d}{ds}\E(g(t),K(s))|_{s=t}=0$
for a.e.\ $t\in[0,1]$. In the special case $g(t)=t\,h$ for a given
function $h\in H^1(\Om)$, we will see (Proposition~\ref{th})
that the last condition
implies  the third condition  considered in Definition~2.9 of
\cite{FraMar3}: $\E(g(t), K(t))\le \E(g(t), K(s))$ for $s<t$.

In Section~\ref{irrev} we prove the following existence result.

\begin{theorem}\label{kt0}
Let $g\in AC([0,1];H^1(\Om))$ and let $K_{0}\in\Kmf$. Then
there exists a function $K\colon [0,1]\to\Kmf$ such that
\smallskip
\begin{itemize}
\item[(a)]
\hfil $\displaystyle \vphantom{\frac{d}{ds}}
K_0\subset K(s)\subset K(t)$ for $0\le s\le t\le 1$, \hfil
\item[(b)] \hfil $\displaystyle \vphantom{\frac{d}{ds}}
\E(g(0),K(0))\leq \E(g(0),K)
\quad\forall \,K\in\Kmf,\,\ K\supset K_0$,\hfil
\item[(c)] \hfil $\displaystyle \vphantom{\frac{d}{ds}}
\hbox{for  }\, 0<t \le1\quad\E(g(t),K(t))\leq \E(g(t),K)
\quad\forall \, K\in\Kmf,\,\ K\supset
{\textstyle\bigcup_{s<t}K(s)}$,\hfil
\item[(d)]\hfil $\displaystyle \vphantom{\frac{d}{ds}}
t\mapsto \E(g(t),K(t)) \hbox{ is
absolutely continuous on }[0,1]$, \hfil
\item[(e)]\hfil$\displaystyle
\frac{d}{ds}\E(g(t),K(s))\Big|\lower1.5ex\hbox{$\scriptstyle s=t$}=0
\quad \hbox{for a.e.\ }t\in[0,1]$.\hfil
\end{itemize}
\smallskip
Moreover every function $K\colon[0,1]\to \Kmf$
which satisfies (a)--(e) satisfies also
\begin{itemize}
\item[(f)] \hfil$\displaystyle\frac{d}{dt}\E(g(t),K(t))=
2\int_{\Om\setminus K(t)} \nabla u(t)\,\nabla \gdot(t)\,dx
\quad \hbox{ for a.e.\ } t\in[0,1]$,\hfil
\end{itemize}
where $u(t)$ is a solution of the minimum problem (\ref{e0}) which defines
$\E(g(t),K(t))$, and $\gdot(t)$ is the time derivative of the function
$g(t)$.
\end{theorem}

If $g(0)=0$, we can prove
that there exists a solution of problem (a)--(e) with $K(0)=K_{0}$
(Remark~\ref{g0}).
We underline that, although we can not exclude that the surface energy
$\huno(K(t))$ may present some jump discontinuities in time
(see \cite[Section 4.3]{FraMar3}), in our result
{\it the total energy is always an
absolutely continuous function of time\/} by condition~(d).

If $\partial_D\Om$ is sufficiently smooth, we can integrate  by parts
the right hand side of (f) and, taking into account the Euler equation
satisfied by $u(t)$, we obtain
\begin{equation}\label{energy}
\frac{d}{dt}\E(g(t),K(t))=
2\int_{\partial_D\Om\setminus K(t)}
\frac{\partial u(t)}{\partial \nu}\, \gdot(t)\,d\huno
\quad \hbox{ for a.e.\ } t\in[0,1]\,,
\end{equation}
where $\nu$ is the outer unit normal to $\partial\Om$. Since the
right hand side of (\ref{energy}) is the power of the force exerted
on the boundary to obtain the displacement $g(t)$ on
$\partial_D\Om\setmeno K(t)$,
equality (\ref{energy}) expresses the {\it conservation of energy\/}
in our quasi-static model, where all kinetic effects are neglected.

The proof of this existence result is obtained by a time 
discretization process.
Given a time step $\delta>0$, for every integer $i\ge 0$ we set
$t_i^\delta:=i\delta$ and $g_i^\delta:=g(t_i^\delta)$.
We define $K_i^\delta$, inductively, as a solution of the minimum problem
\begin{equation}\label{pidelta0}
\min_K\big\{ \E(g_i^\delta,K) : K \in \Kmf,\
K \supset K_{i-1}^\delta \big\}\,,
\end{equation}
where we set $K_{-1}^\delta=K_0$.

Let $u_i^\delta$ be a solution of the minimum problem (\ref{e0})
which defines $\E(g_i^\delta,K_i^\delta)$. On
$[0,1]$ we define  the step
functions $K_\delta$ and $u_\delta$ by setting
$K_\delta(t):=K_{i}^\delta$ and
$u_\delta(t):=u_{i}^\delta$ for $t_i^\delta\le t<t_{i+1}^\delta$.

Using a standard monotonicity argument, we prove that there exists a
sequence $(\delta_k)$ converging to $0$ such that, for every $t\in[0,1]$,
$K_{\delta}(t)$ converges to a compact set $K(t)$ in the Hausdorff
metric as $\delta\to 0$ along this sequence. Then we can apply the results
on the convergence of the solutions to mixed boundary
value problems in cracked domains
established in Section~\ref{convergence}, and we prove that, if $u(t)$ is a
solution of the minimum problem (\ref{e0})
which defines $\E(g(t),K(t))$, then
$\nabla u_{\delta}(t)$ converges to $\nabla u(t)$ strongly in
$L^2(\Om,\R^2)$ as $\delta\to 0$ along the same sequence considered
above.

The technical hypothesis that the sets $K_{\delta}(t)$
have no more than $m$ connected components
plays a crucial role here. Indeed, if this hypothesis is
dropped, the convergence in the Hausdorff
metric of the cracks $K_{\delta}(t)$ to the crack $K(t)$
does not imply the convergence of the corresponding
solutions of the minimum problems, as shown by
many examples in homogenization theory, that can be found, e.g.,
in \cite{Khr},
\cite{Mur},  \cite{Dam},  \cite{Att-Pic}, \cite{Cor}. These
papers show also that this hypothesis would not
be enough in dimension larger than two.

The results of Section~\ref{convergence} are related to those
obtained by A. Chambolle and F. Doveri in \cite{Ch-D} and by D. Bucur
and N. Varchon in \cite{BucVar1} and \cite{BucVar},
which deal with the case of a pure Neumann
boundary condition. Since we impose a Dirichlet boundary condition
on $\partial_D\Om\setmeno K_{\delta}(t)$ and a Neumann
boundary condition on the
rest of the boundary, our results can not be deduced easily from
these papers, so we give an independent proof, which uses the
duality argument of \cite{BucVar}.

{}From this convergence result and from an approximation lemma with
respect to the Hausdorff metric, proved in Section~\ref{Hausdorff}, we
obtain properties (a), (b), (c), (e), and (f) in integrated form,
which implies (d).

The time discretization process described above
turns out to be a useful tool for the proof of the
existence of a solution $K(t)$ of the problem considered in
Theorem~\ref{kt0}, and provides also an effective way for the numerical
approximation of this solution (see \cite{BFM}),
since many algorithms have been
developed for the numerical solution of minimum problems of the form
(\ref{pidelta0}) (see, e.g., \cite{BZ}, \cite{Rich}, \cite{RM}, \cite{Bou},
\cite{Ch}, \cite{Bou-Cha}).

In Section \ref{tips} we study the motion of the tips of the time
dependent
crack $K(t)$ obtained in Theorem~\ref{kt0}, assuming that, in some
open interval $(t_0,t_1)\subset [0,1]$, the crack $K(t)$ has a fixed
number of tips, that these tips move smoothly, and that their paths
are simple, disjoint, and do not intersect $K(t_0)$. We prove
(Theorem~\ref{Griffith}) that in this case
{\it Griffith's criterion for crack growth\/} is valid in our model: the
absolute value of the
{\it  stress intensity factor\/}
(see Theorem~\ref{Grisvard1} and Remark~\ref{stress})
of the solution $u(t)$ is less than or equal
to $1$ at each tip for every $t\in (t_0,t_1)$, and it is equal to $1$
at a given tip for almost every instant $t\in (t_0,t_1)$ in which the 
tip moves with
positive velocity.
\end{section}

\begin{section}{NOTATION AND PRELIMINARIES}\label{notation}

Given an open subset $A$ of $\R^2$, we say that $A$ has a Lipschitz
boundary at a point $x\in\partial A$ if $A$
is the sub-graph of a Lipschitz function near $x$, in the sense that there exist
an orthogonal coordinate system
$(y_1,y_2)$, a rectangle $U=(a,b)\times(c,d)$ containing $x$,
and a Lipschitz function $\Phi\colon(a,b)\to(c,d)$, such that
$A\cap U=\{y\in U: y_2<\Phi(y_1)\}$. The set of all these points $x$
is the {\it Lipschitz part of the boundary\/} and will be denoted by
$\partial_{L}A$.
If $\partial_{L}A=\partial A$, we say that $A$ has a Lipschitz boundary.

Besides the Sobolev space $H^1(A)$ we shall use also the
{\it Deny-Lions space\/}
$\hb(A):=\{u\in L^2_{loc}(A)\;|\; \nabla u\in L^2(A;\R^2)\}$, which
coincides with the space of all distributions $u$
on $A$
such that $\nabla u\in L^2(A;\R^2)$ (see, e.g.,
\cite[Theorem~1.1.2]{Ma}).
For the proof of the following result  we refer, e.g., to
\cite[Section 1.1.13]{Ma}.
\begin{proposition}\label{closed}
The set $\{\nabla u: u\in \hb(A)\}$ is closed in $L^2(A;\R^2)$.
\end{proposition}
Under some regularity assumptions on the boundary, the following
result holds.
\begin{proposition}\label{h1}
Let $u\in \hb(A)$ and $x\in\partial_L A$.
Then there exists a neighbourhood $U$ of $x$ such that
$u|_{A\cap U}\in H^1(A\cap U)$. In particular, if $A$ is bounded and
has  a Lipschitz boundary, then $\hb(A)= H^1(A)$.
\end{proposition}
\begin{proof}
Let $U$ be the rectangle given by the definition of Lipschitz
boundary. It is easy to check that $A\cap U$ has a Lipschitz
boundary. The
conclusion follows now from the Corollary
to Lemma 1.1.11 in \cite{Ma}.
\end{proof}

We recall some properties of the functions in the spaces
$H^1(A)$ and $\hb(A)$, which are related to the notion of capacity.
For more details we refer to \cite{EG}, \cite{HKM}, \cite{Ma}, and \cite{Zie}.

\begin{definition}
Let $B$ be a bounded open set in $\R^2$.
The capacity of an arbitrary subset $E$ of $B$ is defined as
$$
{\rm cap} (E,B) :=\inf_{u\in{\mathcal U}^B_{E}}
\int_{B}|\nabla u|^2\,dx \,,
$$
where ${\mathcal U}^B_{E}$ is the set of all functions $u\in H^1_0(B)$
such that $u\geq 1$ a.e.\ in a neighbourhood of~$E$.
\end{definition}
We say that a property is true {\it quasi-everywhere\/} on a set
$E\subset B$, and
write {\it q.e.\/}, if it holds on $E$ except on a set of capacity 
zero. As usual,
the expression almost everywhere, abbreviated as a.e., refers to the
Lebesgue measure.
A function $u\colon E\to\overline \R$ is said to be {\it
quasi-continuous\/} on $E$ if
for every
$\e>0$ there exists an open set $U_{\e}$, with ${\rm cap}(U_{\e},B)<\e$, such
that $u|_{E\setminus U_{\e}}$ is continuous on $E\setmeno U_{\e}$.
It is easy to prove that both notions of  quasi-everywhere and
quasi-continuity do not depend on~$B$.

It is known that every function $u\in\hb(A)$ has a {\it quasi-continuous
representative\/} $\tilde u$, which is uniquely defined q.e.\ on
$A\cup\partial_L A$, and satisfies
\begin{equation}\label{medie}
\lim_{\rho\to 0} \;
-\hskip-1.1em\int_{B_{\rho}(x)\cap A}|u(y)-\tilde u(x)|\,dy=0\qquad
\hbox{for q.e.\ } x\in A\cup\partial_L A\,,
\end{equation}
where $-\hskip-.9em\int$ denotes the average and $B_\rho(x)$ is the
open ball with centre $x$ and radius $\rho$.
If $u_n\to u$ strongly in $H^1(A)$, then a
subsequence of $(\tilde u_n)$ converges to $\tilde u$ q.e. in
$A\cup\partial_L A$.
If $u,\, v\in \hb(A)$ and their traces coincide $\huno$-a.e.\ on
$\partial_L A$, then $\tilde u$ and $\tilde v$
coincide q.e.\ on $\partial_L A$.

In the quoted books the quasi-continuous representatives are defined
only on $A$. The straightforward definition of $\tilde u$ on $\partial_L A$
relies on the existence of extension operators for Lipschitz domains;
the q.e.\ uniqueness of $\tilde u$ on $\partial_L A$ can be deduced from
(\ref{medie}).
To simplify the notation we shall always identify each function $u\in
\hb(A)$ with its quasi-continuous representative $\tilde u$.

Propositions~\ref{closed} and \ref{h1} imply the following result.
\begin{corollary}\label{complete} Assume that $A$ is connected, and
let ${\Gamma}$ be a non-empty relatively open subset of $\partial A$
with ${\Gamma}\subset\partial_L A$. Then the space
$\hb_0(A,\Gamma):=\{u\in\hb(A):u=0\;\hbox{ q.e.\ on }{\Gamma}\}$ is a 
Hilbert space
with the norm $\|\nabla u\|_{L^2(A;\R^2)}$. Moreover, if $(u_n)$ is a
bounded sequence in $\hb_0(A,\Gamma)$, then there exist a subsequence,
still denoted by $(u_n)$, and a function $u\in \hb_0(A,\Gamma)$ such
that $\nabla u_n\wto \nabla u$ weakly in $L^2(A;\R^2)$.
\end{corollary}
\begin{proof}
Let $(v_n)$ be  a Cauchy sequence in $\hb_0(A,\Gamma)$.
We can construct an increasing sequence
$(A_k)$ of connected open sets with
Lipschitz boundary such that
$A=\bigcup_k A_k$, and
${\Gamma}=\bigcup_{k}(\partial A_k\cap\partial A)$.

By Proposition~\ref{h1} the functions $v_n$ belong to
$H^1(A_k)$ and $v_n=0$ q.e.\ on
$\partial A_k\cap\partial A$. As
$\huno(\partial A_k\cap\partial A)>0$ for $k$ large enough,
by the Poincar\'e inequality
$(v_n)$ is a Cauchy sequence in $H^1(A_k)$, and therefore it converges
  strongly in $H^1(A_k)$ to a function $v$ with $v=0$ q.e.\ on
  $\partial A_k\cap\partial A$.  It is then easy to construct a function
  $v\in \hb(A)$ such that $v=0$ q.e.\ on ${\Gamma}$ and
$v_n\to v$ strongly in $H^1(A_k)$ for every $k$. As
  $(\nabla v_n)$ converges strongly in
$L^2(A;\R^2)$, we conclude that
$v_n\to v$ strongly in
$\hb_0(A;\Gamma)$.

Let $(u_n)$ be a
bounded sequence in $\hb_0(A,\Gamma)$. As in the previous part of the
proof we deduce that $(u_n)$ is bounded in $H^1(A_k)$ for every $k$.
By a diagonal argument we can prove that there exist a subsequence,
still denoted by $(u_n)$, and a function $u\in \hb(A)$ such
that $u_n\wto u$ weakly in $H^1(A_k)$ for every $k$. Then a sequence 
of convex combinations of the functions $u_n$ converges to $u$ 
strongly in $H^1(A_k)$. This implies
$u=0$ q.e.\ on
  $\partial A_k\cap\partial A$ for every $k$, hence
  $u\in \hb_0(A,\Gamma)$. As
  $(\nabla u_n)$ is bounded in
$L^2(A;\R^2)$, we conclude that
$\nabla u_n\wto \nabla u$
weakly in $L^2(A;\R^2)$.
\end{proof}

\begin{proposition}\label{dueconnessi}
Let $u\in \hb(A)$ and let $C_1$ and $C_2$ be two connected
subsets of $A\cup\partial_LA$
with $\overline C_1\cap \overline C_2\neq \emptyset$.
Assume that $u$ is constant q.e.\ on $C_i$  for $i=1,\,2$. Then $u$ is
constant q.e.\ on $C_1\cup C_2$.
\end{proposition}

\begin{proof} We may assume that $C_1$ and $C_2$ have more than one
point, since otherwise the statement is trivial.
Let us denote the constant values of $u$ on $C_1$ and
$C_2$ by $c_1$ and $c_2$ respectively, and
let us fix  $x\in \overline C_1\cap \overline C_2$.
Since $x\in A\cup\partial_LA$, we may assume that $u$ belongs to
$H^1(B_r(x))$ for some $r>0$ (we use an extension operator if
$x\in \partial_LA$), and that $C_i\cap\partial B_\rho(x)\neq \emptyset$
for  $i=1,\,2$ and
$0<\rho<r$.
Hence for almost every $\rho\in (0,r)$ (the quasi-continuous
representative of) $u$ takes the values
$c_1$ and $c_2$ in two distinct points of $\partial B_\rho(x)$. This
implies
$$
\int_{\partial B_\rho(x)}|\nabla u|^2\,d\huno\ge
\frac{(c_2-c_1)^2}{\pi\rho}\,,
$$
which yields $\nabla u\notin L^2(B_r(x);\R^2)$,
in contradiction with
our assumption.
\end{proof}

We conclude this section by stating a property of connected sets with 
finite length.

\begin{proposition}\label{connected}
Let $C$ be a connected subset of $\R^2$. Then $\huno(\overline C)=\huno(C)$.
\end{proposition}

\begin{proof}
It is clearly enough to prove the statement when $\huno(C)<+\infty$.
The following concise argument was suggested by Luigi Ambrosio.
If $x$, $y\in C$, then 
$\huno(C)\geq |x-y|$ (the classical proof, see e.g., \cite[Lemma 
3.4]{Fal}, does not need the hypothesis that $C$ is compact). 
Therefore $\huno(C\cap B_\rho(x))\ge \rho$ for every $x\in \overline C$ 
and $0<\rho<{\rm diam}(C)/2$. This 
implies that $\huno({\overline C\setmeno C})=0$
by a standard argument based on the Besicovitch 
covering lemma (see \cite[2.10.19(4)]{Fed}).
\end{proof}

\end{section}

\begin{section}{HAUSDORFF MEASURE AND HAUSDORFF CONVERGENCE}
\label{Hausdorff}

Throughout the paper $\Om$ is a fixed {\it bounded connected open\/}
subset of $\R^{2}$ with {\it Lipschitz boundary\/}.
In this section we study the
behaviour of the Hausdorff measure
$\huno$ along suitable sequences of compact sets which converge in the
Hausdorff metric.

Let $\K$ be the set of all compact subsets
of $\overline\Om$, and let
$\Kf:=\{{K\in \K}: {\huno(K)<+\infty}\}$. Given an integer $m\ge 1$, let
$\Km$ be the set of all compact subsets of $\overline\Om$
with at most $m$ connected components, and let
$\Kmf:=\{{K\in \Km}: {\huno(K)<+\infty}\}$. For every $\lambda\ge0$
we consider also the set $\Kml:=\{{K\in \Km}: {\huno(K)\le \lambda}\}$.

We recall that the {\it Hausdorff distance\/} between
$K_{1},\, K_{2}\in\K$ is defined by
$$
d_{H}(K_{1},K_{2}):=
\max\big\{ \sup_{x\in K_1}{\rm dist}(x,K_2),
\sup_{y\in K_2}{\rm dist}(y,K_1)\big\}\,,
$$
with the conventions  ${\rm dist}(x,\emptyset)={\rm diam}(\Om)$ and
$\sup\emptyset=0$, so that $d_{H}(\emptyset, K)=0$ if $K=\emptyset$
and $d_{H}(\emptyset, K)={\rm diam}(\Om)$ if $K\neq\emptyset$.
We say that $K_n\to K$ in the Hausdorff metric if
${d_H(K_n,K)\to0}$. The following compactness theorem is
well-known (see, e.g., 
\cite[Blaschke's Selection Theorem]{Rog}).
\begin{theorem}\label{compactness}
Let $(K_n)$ be a sequence in $\K$.
Then there exists a subsequence which
converges in the Hausdorff metric to a set
$K\in \K$.
\end{theorem}


It is well-known that, in general, the Hausdorff measure is not
lower semicontinuous on $\K$
with respect to the convergence in the
Hausdorff metric. When all sets are connected, we have the
following lower semicontinuity theorem, whose proof can be obtained
as in Theorem 10.19 of \cite{MS}.

\begin{theorem}[Go\l \c ab's Theorem]\label{Golab}
Let $(K_n)$ be a sequence in
$\Kone$ which converges to $K$ in the Hausdorff metric. Then
$K\in \Kone$ and
$$
\huno(K\cap U)\le \liminf_{n\to\infty} \,\huno(K_n\cap U)
$$
for every open set $U\subset \R^2$.
\end{theorem}

Go\l \c ab's Theorem says that, for every $\lambda<+\infty$, $\Konel$
is closed under convergence in the Hausdorff metric. In the next
corollary we extend this result to $\Kml$.
\begin{corollary}\label{Golab2}
Let $m\ge 1$ and let $(K_n)$ be a sequence  in
$\Km$ which converges to $K$ in the Hausdorff metric. Then
$K\in \Km$ and
$$
\huno(K\cap U)\le \liminf_{n\to\infty} \,\huno(K_n\cap U)
$$
for every open set $U\subset \R^2$.
\end{corollary}

\begin{proof}
Let $K^{1}_n,\ldots, K^{k_n}_n$ be the connected components of
$K_n$. As $k_n\le m$, there exists $k\le m$ such that, up to a
subsequence, $k_n=k$ for all~$n$. By Theorem \ref{compactness}
we may also assume that
$K^{1}_n\to \widehat K{}^{1}$, $\ldots$, $K^{k}_n\to \widehat K{}^{k}$
in the Hausdorff  metric, where $\widehat K{}^{1},\ldots, \widehat K{}^{k}$
are compact and connected.

We claim that
\begin{equation}\label{claim}
K\subset \widehat K{}^{1}\cup \cdots \cup \widehat K{}^{k}\,.
\end{equation}
Indeed,
for every $x\in K$ there exists a sequence $x_n\to x$ such that
$x_n\in K_n$, which implies $x_n\in K^{i_n}_n$ for some $i_n$ between $1$
and $k$. Hence there exists $i$ such that $i_n=i$ for infinitely many
indices $n$, and, consequently, $x\in \widehat K{}^{i}$.
This proves (\ref{claim}),
which implies that $K$ has at most $k\le m$ connected components.

By Go\l \c ab's Theorem \ref{Golab} we have
$$
\huno(\widehat K^{j}\cap U)\le
\liminf_{n\to\infty} \,\huno(K^{j}_n\cap U)
$$
for $j=1,\ldots,k$. The conclusion follows now from (\ref{claim}).
\end{proof}

We shall use also the following consequence of
Corollary~\ref{Golab2}.
\begin{corollary}\label{sci2}
Let $(H_n)$ be a sequence in $\K$ which
converges to $H$ in the Hausdorff metric.
Let $m\ge 1$ and let
$(K_n)$ be a sequence in $\Km$  which
converges to $K$ in the Hausdorff metric.
Then
\begin{equation}\label{scidif}
\huno(K\setmeno H)\leq
\liminf_{n\to\infty} \, \huno(K_n\setmeno H_n)\,.
\end{equation}
\end{corollary}
\begin{proof}
Given $\e>0$, let
$H^{\e}:=\{x\in\overline\Om: {\rm dist}(x,H)\le \e\}$. As
$H_n\subset  H^{\e}$ for $n$ large enough, we have
$K_n\setmeno H^{\e}\subset  K_n\setmeno H_n$.
Applying Corollary \ref{Golab2} with $U=\R^2\setmeno H^{\e}$ we get
$$
\huno(K\setmeno H^{\e}) \le
\liminf_{n\to\infty} \, \huno(K_n\setmeno H^{\e})
\le \liminf_{n\to\infty} \, \huno(K_n\setmeno H_n)\,.
$$
Passing to the limit as $\e\to0$ we obtain (\ref{scidif}).
\end{proof}

In Section~\ref{irrev} we shall use the following approximation result.

\begin{lemma}\label{differ}
Let $p$ and $m$ be positive integers, let $(H_n)$ be a sequence in $\Kpf$
which converges in the
Hausdorff metric to $H\in\Kpf$,
and let $K$ be an element of $\Kmf$ with $K\supset H$.
Then there exists a sequence $(K_n)$ in $\Kmf$ such that
$K_n\to K$ in the Hausdorff metric,
$H_n\subset K_n$, and
$\huno(K_n\setmeno H_n)\to\huno(K\setmeno H)$.
\end{lemma}

To prove Lemma \ref{differ} we need the following three lemmas.

\begin{lemma}\label{differ1}
Let $H\in\Kone$ and let $(H_n)$ be a sequence in $\Kp$ which
converges to $H$ in the Hausdorff metric. Then there exists a
sequence $(\widehat H_n)$ in $\Kone$ such that
$\widehat H_n\to H$ in the Hausdorff metric,
$H_n\subset \widehat H_n$, and
$\huno(\widehat H_n\setmeno H_n)\to0$.
\end{lemma}

\begin{proof}
Passing to a subsequence, we may assume,
as in the first part of the proof of Corollary \ref{Golab2}, that there
exists a constant $k\le p$ such that every
$H_n$ has exactly $k$ connected components $H^{1}_n,\ldots, H^{k}_n$
and $H^{1}_n\to\widehat H{}^{1}$, $\ldots$,
$H^{k}_n\to\widehat H{}^{k}$
in the Hausdorff  metric, where $\widehat H{}^{1},\ldots, \widehat H{}^{k}$
are compact and connected and
$H=\widehat H{}^{1}\cup \cdots \cup \widehat H{}^{k}$.

As $H$ is connected, there exists a finite family of indices
$(\sigma_j)_{0\le j\le \ell}$, with
$\{\sigma_0,\ldots,\sigma_\ell\}=\{1,\ldots,k\}$, such that
$\widehat H{}^{\sigma_{j-1}} \cap \widehat H{}^{\sigma_{j}}\neq
\emptyset$ for $j=1,\ldots, \ell$. Let us fix a point $x^j\in
\widehat H{}^{\sigma_{j-1}} \cap \widehat H{}^{\sigma_{j}}$. By
the convergence in the Hausdorff metric there exist
$x^j_n\in H^{\sigma_{j-1}}_n$ and  $y^j_n\in H^{\sigma_{j}}_n$
such that $x^j_n\to x^j$ and $y^j_n\to x^j$ as $n\to\infty$.

Since $\Om$ has a Lipschitz boundary, there exist arcs $X^j_n$ and
$Y^j_n$ in $\overline\Om$, connecting $x^j_n$ to $x^j$ and $y^j_n$ to $x^j$
respectively, such that $\huno(X^j_n)\to 0$ and $\huno(Y^j_n)\to 0$
as $n\to\infty$. Let us define
$$
\widehat H_n:= H_n \cup \bigcup_{j=1}^\ell X^j_n
\cup \bigcup_{j=1}^\ell Y^j_n\,.
$$
It is clear that $\widehat H_n\to H$ in the Hausdorff metric and that
$\huno(\widehat H_n\setmeno H_n)\to0$. Since
$$
\widehat H_n=H^{\sigma_{0}}_n \cup X^1_n \cup Y^1_n \cup
H^{\sigma_{1}}_n\cup \cdots \cup
H^{\sigma_{\ell-1}}_n \cup X^\ell_n \cup Y^\ell_n \cup
H^{\sigma_{\ell}}_n\,,
$$
we conclude that $\widehat H_n$ is connected.
\end{proof}

\begin{lemma}\label{differ2}
Let $K\in\Konef$ and let $H$ be a non-empty compact subset of $K$ with
$p\ge2$ connected components $H^{1},\ldots, H^{p}$.
Then there exist a finite family of indices
$(\sigma_j)_{0\le j\le \ell}$, with
$\{\sigma_0,\ldots,\sigma_\ell\}=\{1,\ldots,p\}$, and a family
$(\Gamma_j)_{1\le j\le \ell}$ of connected components of ${K\setmeno H}$,
such that $\overline \Gamma_j$ connects 
$H^{\sigma_{j-1}}$ with $H^{\sigma_{j}}$ for
$j=1,\ldots,\ell$.
\end{lemma}

\begin{proof}
It is clear that $K\setmeno H\neq\emptyset$, since otherwise $H$ has 
exactly one connected component.
Since $K$ is locally connected (see, e.g.,
\cite[Lemma 1]{Ch-D}), and
${K\setmeno H}$ is open in $K$, the connected components
$C$ of ${K\setmeno H}$ are open in $K$.
Since each $C$ is closed in $K\setmeno H$, we have
$C=\overline C\cap (K\setmeno H)$.
If $C=\overline C$, then $K$ would contain an open, closed, and
non-empty proper subset (recall that $H\neq\emptyset$),
which contradicts the fact that $K$ is connected.
Therefore $C\neq \overline C$. As
$\overline C\cap ({K\setmeno H})=C$,
we conclude that $\emptyset\neq{\overline C\setmeno C}\subset H$.
Therefore $\overline C\cap H\neq\emptyset$
for every connected component $C$ of ${K\setmeno H}$.

For $j=1,\ldots,p$ let $\widehat K{}^{j}$ be the union of $H^{j}$
and of all the connected components $C$ of ${K\setmeno H}$
such that $\overline C\cap H^{j}\neq\emptyset$. To prove that
$\widehat K{}^{j}$ is open in $K$, we fix a sequence $(x_n)$ in
${K\setmeno \widehat K{}^{j}}$ which converges to a point $x\in K$.
If $x_n\in {H\setmeno H^{j}}$ for infinitely many indices $n$,
then $x\in {H\setmeno H^{j}}$, hence
$x\notin \widehat K{}^{j}$. If there exists a connected component
$C_0$ of
${K\setmeno H}$ such that $\overline C_0\cap H^{j}=\emptyset$ and
$x_n\in \overline C_0$ for infinitely many indices $n$, then
$x\in \overline C_0$; this implies $x\notin \widehat K{}^{j}$, since
$\overline C_0\cap C=\emptyset$ for every connected component
$C$ of ${K\setmeno H}$
with $\overline C\cap H^{j}\neq\emptyset$. In the other
cases there exists a sequence $(C_n)$ of pairwise disjoint
connected components of ${K\setmeno H}$, with
$\overline C_n\cap H^{j}=\emptyset$,
such that, up to a subsequence, $x_n\in \overline C_n$. As
$\huno(K)<+\infty$, by Proposition \ref{connected}
$\huno(\overline C_n)=\huno(C_n)\to 0$, hence ${\rm
dist}(x_n,{H\setmeno H^{j}})\to 0$, which gives $x\in H\setmeno H^{j}$,
so that $x\notin \widehat K{}^{j}$ also in this case.
Therefore $\widehat K{}^{j}$ is open in $K$.

Since $\overline C\cap H\neq\emptyset$
for every connected component $C$ of ${K\setmeno H}$, we have
$K=\widehat K{}^{1}\cup\cdots \cup \widehat K{}^{p}$.
As $K$ is connected, there exists a finite family of indices
$(\sigma_j)_{0\le j\le \ell}$, with
$\{\sigma_0,\ldots,\sigma_\ell\}=\{1,\ldots,p\}$, such that
$\sigma_{j-1}\neq \sigma_{j}$ and
$\widehat K{}^{\sigma_{j-1}} \cap \widehat K{}^{\sigma_{j}}\neq
\emptyset$ for $j=1,\ldots, \ell$. As
$H^{\sigma_{j-1}} \cap H^{\sigma_{j}}=\emptyset$, there exists a
connected component $\Gamma_j$ of ${K\setmeno H}$ such that
$\Gamma_j\subset
\widehat K{}^{\sigma_{j-1}} \cap \widehat K{}^{\sigma_{j}}$
and, consequently, $H^{\sigma_{j-1}}\cap \overline \Gamma_j\neq\emptyset\neq
H^{\sigma_{j}}\cap \overline \Gamma_j$.
\end{proof}

\begin{lemma}\label{differ3}
Let $p$ be a positive integer, let $(H_n)$ be a sequence in $\Kpf$
which converges in the
Hausdorff metric to $H\in\Kpf$,
and let $K$ be an element of $\Konef$ with $K\supset H$.
Then there exists a sequence $(K_n)$ in $\Konef$ such that
$K_n\to K$ in the Hausdorff metric,
$H_n\subset K_n$, and
$\huno(K_n\setmeno H_n)\to\huno(K\setmeno H)$.
\end{lemma}

\begin{proof}
If $H=\emptyset$, we just define $K_n:=K$ and notice that
$H_n=\emptyset$ for $n$ large enough.

Assume now $H\neq\emptyset$ and let $H^{1},\ldots, H^{k}$, $k\le
p$, be its connected components. If $k=1$ we set
$\widehat K := H=H^{1}$. If $k\ge 2$, by Lemma \ref{differ2}
there exist a finite family of indices
$(\sigma_j)_{0\le j\le \ell}$, with
$\{\sigma_0,\ldots,\sigma_\ell\}=\{1,\ldots,k\}$, and a family
$(\Gamma_j)_{1\le j\le \ell}$ of connected components of ${K\setmeno H}$,
such that $H^{\sigma_{j-1}}\cap \overline \Gamma_j\neq\emptyset\neq
H^{\sigma_{j}}\cap \overline \Gamma_j$ for $j=1,\ldots,\ell$; in
this case we set
$$
\widehat K := H \cup\bigcup_{j=1}^\ell \overline\Gamma_j\,.
$$

In both cases we want to construct a sequence $(\widehat K_n)$ in
$\Konef$ which converges to $\widehat K$ in the Hausdorff metric and
such that $H_n\subset\widehat K_n$ and
\begin{equation}\label{khat}
\limsup_{n\to\infty} \huno(\widehat K_n\setmeno H_n)\le
\huno(\widehat K\setmeno H)\,.
\end{equation}

Let us fix $\e>0$ such that the sets
$\{{x\in\overline\Om}: {{\rm dist}(x, H^{i})\le \e}\}$,
$i=1,\ldots,k$, are pairwise disjoint, and let
$$
\widetilde H{}^{i}_n:=
\{{x\in H_n}: {{\rm dist}(x, H^{i})\le \e}\}\,.
$$
It is easy to see that $\widetilde H{}^{i}_n\in\Kpf$ and
$H_n=	\widetilde H{}^{1}_n \cup\cdots\cup \widetilde H{}^{k}_n$
for $n$ large enough, and that $(\widetilde H{}^{i}_n)$ converges
to $H^{i}$ in the Hausdorff metric  as $n\to \infty$. By Lemma
\ref{differ1} there exists a
sequence $(\widehat H{}^{i}_n)$ in $\Konef$ such that
$\widehat H{}^{i}_n\to H^{i}$ in the Hausdorff metric,
$\widetilde H{}^{i}_n\subset \widehat H{}^{i}_n$, and
$\huno(\widehat H{}^{i}_n\setmeno \widetilde H{}^{i}_n)\to0$.

If $k=1$ we define $\widehat K_n:= \widehat H{}^{1}_n$.

If $k\ge 2$, for every $j=1,\ldots,\ell$
we fix two points $x^j\in H^{\sigma_{j-1}} \cap \overline\Gamma_j$
and $y^j\in H^{\sigma_{j}} \cap \overline\Gamma_j$. By
the convergence in the Hausdorff metric there exist
$x^j_n\in \widehat H{}^{\sigma_{j-1}}_n$  and
$y^j_n\in \widehat H{}^{\sigma_{j}}_n$
such that $x^j_n\to x^j$ and $y^j_n\to y^j$ as $n\to\infty$.
Since $\Om$ has a Lipschitz boundary, there exist arcs $X^j_n$ and
$Y^j_n$
in $\overline\Om$, connecting $x^j_n$ to $x^j$ and
$y^j_n$ to $y^j$ respectively, such that $\huno(X^j_n)\to 0$ and
$\huno(Y^j_n)\to 0$
as $n\to\infty$.
Let us define
$$
\widehat K_n:= \bigcup_{i=1}^k \widehat H{}^{i}_n
\cup \bigcup_{j=1}^\ell X^j_n
\cup \bigcup_{j=1}^\ell \overline\Gamma_j
\cup \bigcup_{j=1}^\ell Y^j_n\,.
$$

In both cases $k=1$ and $k\ge 2$
it is clear that $\widehat K_n\to \widehat K$ in the Hausdorff metric
and that (\ref{khat}) holds, since by Proposition~\ref{connected}
$\huno(\overline\Gamma_j)=\huno(\Gamma_j)$. As
$\widehat K_n=\widehat  H{}^{1}_n$ for $k=1$, and
$$
\widehat K_n=\widehat  H{}^{\sigma_{0}}_n \cup X^1_n
\cup \overline\Gamma_1 \cup Y^1_n \cup
\widehat H{}^{\sigma_{1}}_n\cup \cdots \cup
\widehat H{}^{\sigma_{\ell-1}}_n \cup X^\ell_n
\cup \overline\Gamma_\ell \cup Y^\ell_n \cup
\widehat H{}^{\sigma_{\ell}}_n
$$
for $k\ge 2$, we conclude that $\widehat K_n$ is connected in both cases.

As the connected components $C$ of ${K\setmeno \widehat K}$ are
connected components of ${K\setmeno H}$, the argument given
at the beginning of the proof of Lemma \ref{differ2}
shows that each $C$ is open in $K$ and satisfies $\overline C\cap
H\neq\emptyset$. Since $K$ is separable, the connected components
of ${K\setmeno \widehat K}$ form a finite or countable sequence $(C_i)$.

For every $i$ we fix a
point $z^i\in \overline C_i\cap H$.
As $H_n\to H$ in the Hausdorff metric, there exists
$z^i_n\in H_n$ such that
$z^i_n\to z^i$ as $n\to\infty$.
Since $\Om$ has a Lipschitz boundary, for every $i$
there exists an arc $Z^i_n$
in $\overline\Om$, connecting $z^i_n$ to $z^i$,
such that $\huno(Z^i_n)\to 0$
as $n\to\infty$.

If there are infinitely many connected components $C_i$,
there exists a sequence of integers $(h_n)$ tending to
$\infty$ such that
\begin{equation}\label{aax}
\lim_{n\to\infty} \sum_{i=1}^{h_n} \huno(Z^i_n)= 0\,.
\end{equation}
If there are $h<+\infty$ connected components $C_i$, (\ref{aax}) is true with
$h_n=h$ for every $n$.
Let
$$
K_n:=\widehat K_n
\cup \bigcup_{i=1}^{h_n}  Z^i_n
\cup \bigcup_{i=1}^{h_n} \overline C_i
\,.
$$
Then the sets $K_n$ are connected, contain $H_n$,
and converge to $K$ in the
Hausdorff metric. As
$\huno(C_i)=\huno(\overline C_i)$ by Proposition~\ref{connected},
we have
$$
\huno(K_n\setmeno H_n)\le \huno(\widehat K_n\setmeno H_n)+
\sum_{i=1}^{h_n} \huno(Z^i_n)+
\sum_{i=1}^{h_n} \huno(C_i)
\,,
$$
which, together with (\ref{khat}) and (\ref{aax}), yields
\begin{equation}\label{aay}
\limsup_{n\to\infty} \,\huno(K_n\setmeno H_n)\le
\textstyle\huno(\widehat K\setmeno H)+
\huno(\bigcup_{i}C_i)=
\huno(K\setmeno H)\,.
\end{equation}
The opposite inequality for the lower limit follows from
  Corollary~\ref{sci2}.
  \end{proof}

\begin{proof}[Proof of Lemma~\ref{differ}.]
Let $K^{1},\ldots, K^{k}$, $k\le m$, be the connected components
of $K$. Let us fix $\e>0$ such that the sets
$\{{x\in\overline\Om}: {{\rm dist}(x, K^{i})\le \e}\}$,
$i=1,\ldots,k$, are pairwise disjoint, and let
$$
\widehat H{}^{i}_n:=
\{{x\in H_n}: {{\rm dist}(x, K^{i})\le \e}\}\,.
$$
It is easy to see that $\widehat H{}^{i}_n\in\Kpf$ and
$H_n=\widehat H{}^{1}_n \cup\cdots\cup \widehat H{}^{k}_n$
for $n$ large enough, and that $(\widehat H{}^{i}_n)$ converges
to $H^i:=H\cap K^{i}$ in the Hausdorff metric  as $n\to \infty$. By Lemma
\ref{differ3} there exists a
sequence $( K^{i}_n)$ in $\Konef$ such that
$ K^{i}_n\to K^{i}$ in the Hausdorff metric,
$\widehat H{}^{i}_n\subset K^{i}_n$, and
$\huno(K^{i}_n\setmeno \widehat H{}^{i}_n)\to
\huno(K^{i}\setmeno H^i)$.
It suffices now to take $K_n:= K{}^{1}_n\cup\cdots
\cup K{}^{k}_n$.
  \end{proof}

\end{section}

\begin{section}{PROPERTIES OF THE HARMONIC CONJUGATE}\label{conjugate}

In the rest of the paper $\partial_N\Om$ is a fixed (possibly empty)
relatively open subset of $\partial\Om$, with a finite
number of connected components, on which we impose a Neumann
boundary condition. Let
$\partial_D\Om:=\partial\Om\setmeno \overline{\partial_N\Om}$,
which turns out to be a relatively open subset of $\partial\Om$, with a finite
number of connected components. On this set
we want to impose a Dirichlet boundary condition.

Given  $K\in\K$,
we consider the following boundary value problem:
\begin{equation}\label{*}
\left\{\begin{array}{ll}
\Delta u=0 & \hbox{in }\Omk\,,\\
\noalign{\vskip5pt}
\frac{\partial u}{\partial\nu}=0 & \hbox{on }
\partial(\Omk)\cap(K\cup \partial_N\Om)\,.
\end{array}\right.
\end{equation}
By a solution of (\ref{*}) we mean a function $u$ which satisfies the
following conditions:
\begin{equation}\label{**}
\left\{\begin{array}{l}
u\in \hb(\Omk)\,,\\
\noalign{\vskip5pt}
\displaystyle\int_{\Om\setminus K}\nabla u\,\nabla z\,dx=0\quad\forall z\in
\hb(\Omk)\,,\ z=0\quad \hbox{q.e.\ on } \partial_D\Om\setmeno K\,.
\end{array}\right.
\end{equation}

Since no boundary condition is prescribed on $\partial_D\Om\setmeno 
K$, we do not
expect a unique solution to problem (\ref{*}). Given $g\in \hb(\Omk)$,
we can prescribe the Dirichlet boundary condition
\begin{equation}\label{**b}
u=g\quad\hbox{q.e.\ on }\partial_D\Om\setmeno K\,.
\end{equation}
It is clear that problem (\ref{**}) with the boundary condition
(\ref{**b}) can be solved separately in
each connected component of $\Omk$. By Corollary~\ref{complete} and by
the Lax-Milgram lemma there exists
a unique solution in those components whose boundary meets
$\partial_D\Om\setmeno K$, while on the other components the solution is
given by an arbitrary constant. Thus the solution is not unique, if
there is a connected component whose boundary does not meet
$\partial_D\Om\setmeno K$. Note, however, that $\nabla u$ is always unique.
Moreover, the map $g\mapsto\nabla u$ is linear from $\hb(\Omk)$
into $L^2(\Omk;\R^2)$ and satisfies the estimate
$$
\int_{\Om\setminus K}|\nabla u|^2\,dx\le
\int_{\Om\setminus K}|\nabla g|^2\,dx\,.
$$

By standard arguments on the minimization of quadratic forms it is
easy to see that $u$ is a solution of problem (\ref{**}) and
satisfies the boundary condition (\ref{**b}) if and only if $u$ is a
solution of the minimum problem
\begin{equation}\label{minpb}
\min_{v\in {\mathcal V}(g,K)}\int_{\Om\setminus K}|\nabla v|^2\,dx\,,
\end{equation}
where
\begin{equation}\label{vgk}
{\mathcal V}(g,K):=\{v\in\hb(\Om\setmeno K): v=g\quad \hbox{q.e.\ on }
\partial_D\Om\setmeno K\}\,.
\end{equation}
Throughout the paper, given a function $u\in \hb(\Omk)$ for some 
$K\in\K$, we always
extend $\nabla u$ to $\Om$ by setting $\nabla u=0$ a.e.\ on $K$. Note that,
however, $\nabla u$ is the  distributional gradient  of $u$ 
only in $\Omk$, and, in
general, it does not coincide in $\Om$ with the gradient
of an extension of $u$.

To study the continuous dependence on $K$ of the solutions of problem
(\ref{**}) with boundary condition (\ref{**b}), we shall use the
following lemma.

\begin{lemma}\label{glinf}
Let $(K_n)$ be a sequence  in
$\K$ which converges to $K$ in the Hausdorff
metric.
Let $u_n\in \hb (\Om\setmeno K_n)$ be a sequence such that $u_n=0$
q.e.\ on ${\partial_D\Om\setmeno K_n}$ and
$(\nabla u_n)$ is bounded in
$L^2(\Om;\R^2)$. Then there
exist a subsequence, still denoted by $(u_n)$, and a function
$u\in \hb(\Omk)$, such that $u=0$
q.e.\ on ${\partial_D\Om\setmeno K}$ and $\nabla u_n\wto \nabla u$
weakly in $L^2(U;\R^2)$ for every open set $U\subset\subset\Omk$.
If, in addition, ${\rm meas}(K_n)\to{\rm meas}(K)$, then
$\nabla u_n\wto \nabla u$ weakly in $L^2(\Om;\R^2)$.
\end{lemma}
\begin{proof} Let $C$ be a connected component of $\Omk$ and
let $x\in C$.
Given $0<\e<{\rm dist}(x,\partial C)$, let
$N^{\e}:=\{{x\in \R^2}: {\rm dist}(x,{\partial_N\Om\cup K})\le \e\}$
and let $C^{\e}$ be the connected component of ${C\setmeno N^{\e}}$
containing $x$. For $n$ large enough we have $K_n\subset N^{\e}$.

If the boundary of $C$ meets $\partial_D\Om\setmeno K$,
let ${\Gamma}^{\e}$ be the relative interior of
$\partial C^\e\cap\partial_D\Om$ in $\partial C^\e$.
Since $\partial C$ meets $\partial_D\Om\setmeno K$, for $\e$ small enough
we have ${\Gamma}^{\e}\neq\emptyset$.
As $u_n=0$ q.e.\ on ${\Gamma}^{\e}$,
we apply Corollary~\ref{complete} and deduce that there exists a
function $u\in\hb(C^\e)$, with $u=0$ q.e.\ on ${\Gamma}^{\e}$,
such that, up to a subsequence, $\nabla u_n\wto \nabla u$ weakly in
$L^2(C^\e;\R^2)$.
Since $\e>0$ is arbitrary and $C=\bigcup_\e C_\e$, we can construct
$u\in \hb(C)$, with
$u=0$ q.e.\ on $({\partial C\cap \partial_D\Om})\setmeno K$, such that
$\nabla u_n\wto \nabla u$
weakly in $L^2(U;\R^2)$ for every open set $U\subset\subset C$.

If the boundary of $C$ does not meet $\partial_D\Om\setmeno K$,
passing to a subsequence, we can still assume  that $(\nabla u_n)$
converges weakly in $L^2(C^{\e};\R^2)$ to some function
$\varphi\in L^2(C^{\e};\R^2)$.
  Since the space
$\{\nabla v: v\in \hb (C^{\e})\}$ is closed in $L^2(C^{\e};\R^2)$, we
conclude that there exists $u\in \hb (C^{\e})$ such that
$\nabla u=\varphi$ a.e.\ in $C^{\e}$, and, as in the previous case,
we can construct  $u\in \hb(C)$ such that $\nabla u_n\wto \nabla u$
weakly in $L^2(U;\R^2)$ for every open set $U\subset\subset C$.

Therefore we have constructed
$u\in\hb(\Omk)$, with $u=0$ q.e.\ on ${\partial_D\Om\setmeno K}$,
such that $\nabla u_n\wto \nabla u$
weakly in $L^2(U;\R^2)$ for every open set $U\subset\subset \Omk$.

Assume now that ${\rm meas}(K_n)\to{\rm meas}(K)$ and let $\psi\in
L^2(\Om;\R^2)$. For every $\e>0$ there exists $\delta>0$ such that
$\int_A|\psi|^2\,dx<\e^2$ for ${\rm meas}(A)<\delta$. Let
$U\subset\subset \Omk$ be an open set such that
${\rm meas}((\Om\setmeno K)\setmeno U))<\delta$. As
$U\subset\subset \Om\setmeno K_n$ for $n$ large enough, we have also
${\rm meas}((\Om\setmeno K_n)\setmeno U))<\delta$. Then
$$
\Big|\int_{\Om}(\nabla u_n-\nabla u)\cdot\psi\,dx\Big|\leq
\Big|\int_{U}(\nabla u_n-\nabla u)\cdot\psi\,dx\Big| + c_1\e+c_2\e\,,
$$
where $c_1$ is an upper bound for
$\|\nabla u_n\|_{L^2(\Om;\R^2)}$ and
$c_2:=\|\nabla u\|_{L^2(\Om;\R^2)}$.
{}From the previous part of the lemma
$\limsup_{n}|\int_\Om(\nabla u_n-\nabla u)\cdot\psi\,dx|\leq
c_1\e+c_2\e$ and the conclusion follows from the arbitrariness of~$\e$.
\end{proof}

Throughout the paper $R$ denotes the
rotation on $\R^2$ defined by $R(y_{1},y_{2}):=(-y_{2},y_{1})$.
In the next theorem we prove that every point of $\overline\Om$
has an open neighbourhood $U$ such that every solution $u$ of (\ref{*})
has a harmonic conjugate $v$ on $(U\cap\Om)\setmeno K$
which is constant on each connected component of $\overline U\cap K$ and
on each connected component of $\overline U\cap \partial_N\Om$.

\begin{theorem}\label{conjug} Let $K\in\K$,
let $u$ be a solution of problem (\ref{**}), and let $U$ be an open rectangle
contained in $\Om$, or a rectangle as in
the definition of the Lipschitz part of the boundary. Then
there exists a function $v\in H^1({U\cap\Om})$ such that
$\nabla v=R\,\nabla u$ a.e.\ on ${U\cap\Om}$. Moreover
$v$ is constant q.e.\ on each connected component of
$\overline U\cap K$ and of $\overline U\cap\partial_N\Om$.
\end{theorem}

\begin{proof} If $\varphi\in C^\infty_c(\R^2)$ with ${\rm
supp}(\varphi)\subset {U\cap\Om}$, we have
\begin{equation}\label{div}
\int_{U\cap\Om} \nabla u\,\nabla\varphi\,dx =
\int_{\Om\setminus K} \nabla u\,\nabla\varphi\,dx=0\,,
\end{equation}
where the first equality follows from our convention $\nabla u=0$
a.e.\ in $K$, while the second equality follows from (\ref{**}),
since $\varphi=0$ on $\partial_D\Om$. Equality (\ref{div}) implies that
${\rm div}(\nabla u)=0$ in ${\mathcal D}'({U\cap\Om})$, hence
${\rm rot}(R\,\nabla u)=0$ in ${\mathcal D}'({U\cap\Om})$. As
${U\cap\Om}$ is simply connected and has a Lipschitz boundary,
there exists $v\in H^1({U\cap\Om})$ such that $\nabla v=R\,\nabla u$
a.e.\ in ${U\cap\Om}$.

Since $\nabla v=0$ a.e. in ${U\cap K}$, the function $v$ is constant q.e.\ on
each connected open subset $C$ of ${U\cap K}$, and, by (\ref{medie}),
also on $C\cup\partial_L C$.

To prove that $v$ is constant q.e.\ on
each connected component of $\overline U\cap K$ we use an 
approximation argument.
We write $K$ as the intersection of a decreasing sequence $(K_j)$
of compact subsets of $\overline\Om$ such that
$K\subset {\rm int}_{\overline\Om} K_j$ for every $j$, where
${\rm int}_{\overline\Om} K_j$ denotes the interior of $K_j$
in the relative topology of $\overline\Om$.

Note that $u$ satisfies
\begin{equation}\label{**u}
\left\{\begin{array}{l}
\displaystyle \int_{(U\cap\Om)\setminus K}\nabla u\,\nabla z\,dx=0
\\
\noalign{\vskip5pt}
\forall z\in
\hb((U\cap\Om)\setmeno K)\,,\ z=0\quad\hbox{q.e.\ on }
\partial(U\cap\Om)\setmeno K\,,
\end{array}\right.
\end{equation}
since every such function $z$ can be extended to a function
of $\hb(\Om\setmeno K)$ by setting $z=0$ in
${(\Om\setmeno U)}\setmeno K$.
As $u\in\hb({(U\cap\Om)\setmeno K_j})$, there exists a solution $u_j$ to the
problem
\begin{equation}\label{**j}
\left\{\begin{array}{l}
u_j\in \hb((U\cap\Om)\setmeno K_j)\,,\quad u_j=u\quad\hbox{q.e.\ on }
\partial(U\cap\Om)\setmeno K_j\,,
\\
\noalign{\vskip5pt}
\displaystyle\int_{(U\cap\Om)\setminus K_j}\nabla u_j\,\nabla z\,dx=0
\\
\noalign{\vskip5pt}
\forall z\in
\hb((U\cap\Om)\setmeno K_j)\,,\ z=0\quad\hbox{q.e.\ on }
\partial(U\cap\Om)\setmeno K_j\,.
\end{array}\right.
\end{equation}
Using $u_j-u$ as test function in (\ref{**j}), we obtain that the
norms $\|\nabla u_j\|_{L^2((U\cap\Om)\setminus K_j)}$
are uniformly bounded. By
Lemma~\ref{glinf}, there exists
$u^*\in\hb({(U\cap\Om)\setmeno K})$, with $u^*=u$ q.e.\ on
${\partial(U\cap\Om)\setmeno K}$, such that,  up to a
subsequence,  $(\nabla u_j)$ converges
to $\nabla u^*$ weakly  in $L^2(U\cap\Om;\R^2)$.

Taking $u_j-u^*$ as test function in (\ref{**j}),
we get
$$
\int_{U\cap\Om}|\nabla u_j|^2\,dx=
\int_{U\cap\Om} \nabla u_j\, \nabla u^*\,dx\,.
$$
Passing to the
limit we obtain that $\|\nabla u_j\|_{L^2(U\cap\Om;\R^2)}$ converges to
$\|\nabla u^*\|_{L^2(U\cap\Om;\R^2)}$, hence $\nabla u_j$ converges
to $\nabla u^*$ strongly in
$L^2(U\cap\Om;\R^2)$.

Let us prove that
\begin{equation}\label{uu*}
\nabla u^*=\nabla u\quad\hbox{ a.e.\ in }(U\cap\Om)\setmeno K\,.
\end{equation}
By the uniqueness of the gradients of the solutions,
it is enough to prove that $u^*$ is a solution of (\ref{**u}).

Let $z\in\hb({(U\cap\Om)\setmeno K})$ with $z=0$ q.e.\ on
${\partial(U\cap\Om)\setmeno K}$. As $z\in
\hb({(U\cap\Om)\setmeno K_j})$ and $z=0$ q.e.\ on
${\partial(U\cap\Om)\setmeno K_j}$,
we can use $z$ as test function in (\ref{**j}). Then passing to the limit as
$j\to\infty$ we obtain (\ref{**u}), and the proof of (\ref{uu*}) is
complete.

By the first part of the proof, there exist a function $v_j\in
H^1(U\cap\Om)$, such that $\nabla v_j=R\,\nabla u_j$ a.e.\ on
$U\cap\Om$.
Let $K^{0}$ be a connected component of ${\overline U\cap K}$. It is easy to
see that there exists a connected component $C$ of the interior of
${U\cap K_j}$ such that $K^{0}\subset {C\cup\partial_L C}$ (this is
trivial if $K^{0}\subset{U\cap \Om}$, and follows from the regularity of
$\partial({U\cap \Om})$ if $K^{0}$ meets $\partial({U\cap \Om})$). 
As $v_j$ is constant
q.e.\ on ${C\cup\partial_L C}$, we obtain that $v_j$ is constant
q.e.\ on $K^{0}$.

We may assume that $\int_{U\cap\Om} v_j\,dx=0$ for every $j$.
Since $\nabla v_j=R\,\nabla u_j$ a.e.\ on ${U\cap\Om}$
we deduce that $(\nabla v_j)$ converges to $R\,\nabla u$ strongly in
$L^2({U\cap\Om};\R^2)$, and by the Poincar\'e inequality
$(v_j)$  converges strongly in $H^1({U\cap\Om})$ to a
function $v$ which satisfies
$\nabla v=R\,\nabla u$ a.e.\ on ${U\cap\Om}$. As $v_j$ is constant
q.e.\ on $K^{0}$, we conclude that $v$ is constant
q.e.\ on $K^{0}$.

To prove that $v$ is constant q.e.\ on each connected component of
${\overline U\cap\partial_N\Om}$, it is enough to show that $v$
is constant q.e.\ on  ${V\cap\partial_N\Om}$
whenever $V\subset U$ is a rectangle as in the
definition of the Lipschitz part of the boundary
and ${V\cap\partial\Om}={V\cap\partial_N\Om}$. Let $\psi\in
L^2(V;\R^2)$ be the vector-field defined by $\psi=\nabla u$ a.e.\ in
${V\cap\Om}$ and $\psi=0$ a.e.\ in ${V\setmeno\Om}$.
As at the beginning of the proof, it is easy to see that
${\rm div}(\psi)=0$ in ${\mathcal D}'(V)$, hence
${\rm rot}(R \psi)=0$ in ${\mathcal D}'(V)$. Then
there exists a function $z\in H^1(V)$ such that $\nabla z=R\psi$ a.e\ in $V$.
As $\nabla z=0$ a.e.\ in the connected set ${V\setmeno\overline\Om}$,
using (\ref{medie}) we obtain that $z$ is constant q.e.\ in
${V\setmeno \Om}$.
As $\nabla z=R\,\nabla u=\nabla v$
a.e.\ in the connected set ${V\cap\Om}$, using (\ref{medie}) we obtain
that $z-v$ is constant q.e.\ in ${V\cap\overline\Om}$. {}From these
facts we deduce that $v$ is constant q.e.\ on
$V\cap\partial\Om=V\cap\partial_N\Om$.
\end{proof}

\begin{theorem}\label{due}
Let $K$ be a locally connected compact subset of
$\overline\Om$ and let $u\in \hb({\Om\setmeno K})$.
Assume that for every $x\in\overline\Om$ there exist an open
neighbourhood $U$ of $x$ in $\R^2$ and a function
$v\in H^1({U\cap\Om})$, with
$\nabla v=R\,\nabla u$ a.e.\ in ${U\cap\Om}$,
such that $v$ is constant q.e.\ on
each connected component of $U\cap K$ and of $U\cap\partial_N\Om$.
Then $u$ is a solution of problem~(\ref{**}).
\end{theorem}

\begin{proof}
By a standard localization argument, it is enough to prove that for every
$x\in \overline\Om$ there exists an open neighbourhood $V$ of $x$ in
$\R^2$ such that
\begin{equation}\label{zU}
\left\{\begin{array}{l}
\displaystyle\int_{(V\cap\Om)\setminus K}\nabla u\,\nabla z\,dx=0
\\
\noalign{\vskip5pt}
\forall z\in
\hb((V\cap\Om)\setmeno K)\,,\ z=0\quad \hbox{q.e.\ on }
(V\cap \partial_D\Om)\setmeno K\,,\quad
{\rm supp}(z) \subset\subset V\,.
\end{array}\right.
\end{equation}

For every $x\in \overline\Om$ let $U$ be the neighbourhood
given in the statement of the theorem.
Taking, if necessary, a smaller neighbourhood, we may assume that 
${U\cap\Om}$ has a Lipschitz boundary and that 
$v$ is constant q.e.\ on the closure of
each connected component of ${U\cap K}$ and of ${U\cap\partial_N\Om}$.
Let $V$ be an arbitrary
open neighbourhood of $x$ in
$\R^2$ with $V\subset\subset U$. Since $K$ is locally connected, the
connected components of ${U\cap K}$ are open in $K$, so that only a
finite number of them meets ${\overline V\cap K}$. Similarly, only a
finite number of connected components of ${U\cap \partial_N\Om}$
meets ${\overline V\cap \partial_N\Om}$.
Using Proposition~\ref{dueconnessi} it is easy to prove that there
exist a finite family
$\widehat K{}^{1},\ldots, \widehat K{}^{m}$ of pairwise disjoint
compact sets and a family of distinct constants $c^{1},\ldots, c^{m}$
such that
$$
\overline V\cap (K\cup\partial_N\Om)=
\widehat K{}^{1}\cup \cdots \cup \widehat K{}^{m}
$$
and $v=c^{i}$ q.e.\ on $\widehat K{}^{i}$ for $i=1,\ldots,m$.

We now apply \cite[Theorem~4.5]{HKM}
(to a suitable extension of $v|_{V\cap\Om}$) and
construct a sequence of functions $v_n\in C^{\infty}(\R^2)$, converging
to $v$ in $H^1({V\cap\Om})$, such that $v_n=c^{i}$  in a
neighbourhood $U^{i}_n$ of each $\widehat K{}^{i}$.

Let $z\in \hb({(V\cap\Om)\setmeno K})$, with compact support in $V$,
such that $z=0$ q.e.\ on $(V\cap \partial_D\Om)\setmeno K$, let
$\varphi^{i}_n$ be functions in $C^\infty_c(\R^2)$,
with ${\rm supp}(\varphi^{i}_n)\subset U^{i}_n$, such that
$\varphi^{i}_n=1$ in a neighbourhood of $\widehat K{}^{i}$, and let
$\psi_n:=1-\sum_i \varphi^{i}_n$. By Proposition~\ref{h1} the
function $z\,\psi_n$ belongs to $H^1({(V\cap\Om)\setmeno K})$, and
by \cite[Theorem~4.5]{HKM}  it belongs to
$H^1_0({(V\cap\Om)\setmeno K})$.

Since $\psi_n=1$ where $R\,\nabla v_n\neq 0$, we have
$$
\int_{(V\cap\Om)\setminus K} R\,\nabla v_n\,\nabla z\,dx=
\int_{(V\cap\Om)\setminus K} R\,\nabla v_n\,\nabla (z \,\psi_n)\,dx
=0\,,
$$
where the last equality follows from the fact that
${\rm div}(R\,\nabla v_n)=0$ in $\R^2$ and $z \,\psi_n\in
H^1_0({(V\cap\Om)\setmeno K})$.
Passing to the limit as $n\to\infty$,
we obtain
$$
\int_{(V\cap\Om)\setminus K}\nabla u\,\nabla z\,dx=
-\int_{(V\cap\Om)\setminus K} R\,\nabla v\,\nabla z\,dx=0\,,
$$
showing that $u$ is a solution of (\ref{zU}).
\end{proof}

\end{section}

\begin{section}{CONVERGENCE OF MINIMIZERS}\label{convergence}

In this section we prove the convergence of the minimum points of problems
(\ref{minpb}) corresponding to a sequence $(K_n)$ in $\Kml$ which
converges in the Hausdorff metric.

\begin{theorem}\label{convsol}
Let $m\ge1$ and $\lambda\ge0$, let $(K_n)$ be a sequence
in $\Kml$ which converges to $K$ in the Hausdorff metric, and let
$(g_n)$ be a sequence in $H^1(\Om)$ which converges to
$g$ strongly in $H^1(\Om)$.
Let $u_n$ be a solution
of the minimum problem
\begin{equation}\label{Pn}
\min_{v\in{\mathcal V}(g_n,K_n)}
\int_{\Om\setminus K_{n}}|\nabla v|^2\,dx\,,
\end{equation}
and let $u$ be  a solution of the minimum problem
\begin{equation}\label{P}
\min_{v\in{\mathcal V}(g,K)} \int_{\Om\setminus K}|\nabla v|^2\,dx\,,
\end{equation}
where ${\mathcal V}(g_n,K_n)$ and ${\mathcal V}(g,K)$ are defined by
(\ref{vgk}).
Then $\nabla u_n\to \nabla u$ strongly in $L^2(\Om;\R^2)$.
\end{theorem}

The following lemma is crucial in the proof of Theorem~\ref{convsol}.

\begin{lemma}\label{Mosco}
Let $(K_n)$ be a sequence
in $\Kone$ which converges to $K$ in the Hausdorff metric, and let
$(v_n)$ be a sequence in $H^1(\Om)$ which converges to
$v$ weakly in $H^1(\Om)$. Assume that $v_n=0$ q.e.\ on $K_n$ for
every $n$. Then $v=0$ q.e.\ on $K$.
\end{lemma}

\begin{proof}
Let us fix an open ball $B$ containing $\overline\Om$. Using the same
extension operator we can construct extensions of $v_n$ and $v$,
still denoted by  $v_n$ and $v$, such that $v_n,\,v\in H^1_0(B)$,
$v_n\wto v$ weakly in $H^1(B)$.

Given an open set $A\subset B$, any function $z\in H^1_0(A)$ will
be extended to a function $z\in H^1_0(B)$ by setting $z:=0$ q.e.\ in
$\overline B\setmeno A$. By \cite[Theorem 4.5]{HKM} we have
\begin{equation}\label{h10}
H^1_0(A)=\{z\in H^1(B):z=0\quad\hbox{q.e.\ on }
\overline B\setmeno A\}\,.
\end{equation}
Since the complement of $B\setmeno K_n$ has two connected
components,
from the results of \cite{Sve} and \cite{Buc} we deduce that, for
every $f\in L^2(B)$, the solutions $z_n$ of the Dirichlet problems
$$
z_n\in H^1_0(B\setmeno K_n)\qquad \Delta z_n=f\quad\hbox{in }
B\setmeno K_n
$$
converge strongly in $H^1_0(B)$ to the solution $z$ of the Dirichlet
problem
$$
z\in H^1_0(B\setmeno K)\qquad \Delta z=f\quad\hbox{in }
B\setmeno K\,.
$$
This implies (see, e.g., \cite[Theorem~3.33]{Att})
that, in the space $H^1_0(B)$, the subspaces
$H^1_0(B\setmeno K_n)$ converge to  the subspace
$H^1_0(B\setmeno K)$
in the sense of Mosco (see \cite[Definition 1.1]{Mos}).

Since $v_n\in H^1_0(B\setmeno K_n)$ by (\ref{h10}),
and  $v_n\wto v$ weakly in $H^1(B)$, from the
convergence in the sense of Mosco we deduce that
$v\in H^1_0(B\setmeno K)$, hence $v=0$ q.e.\ on $K$ by (\ref{h10}).
\end{proof}

\begin{proof}[Proof of Theorem~\ref{convsol}.]
Note that $u$ is a minimum point  of (\ref{P}) if and only if $u$ satisfies
(\ref{**}) and (\ref{**b}); analogously, $u_n$ is a minimum point  of
(\ref{Pn})
if and only if  $u_n$ satisfies (\ref{**}) and (\ref{**b}) with $K$
and $g$ replaced
by $K_n$ and $g_n$.

Taking $u_n-g_n$ as test function in the equation satisfied by $u_n$, we
prove that the sequence
$(\nabla u_n)$ is bounded in $L^2(\Om;\R^2)$.
By Lemma~\ref{glinf},
there exists a function $u^*\in \hb(\Omk)$, with $u^*=g$ q.e.\ on
${\partial_D\Om\setmeno K}$,
such that, passing to a subsequence, $\nabla u_n\wto \nabla u^*$
weakly in $L^2(\Om;\R^2)$.

We will prove that
\begin{equation}\label{2uu*}
\nabla u^*=\nabla u\quad\hbox{ a.e.\ in }\Om\setmeno K\,.
\end{equation}
As the limit does not depend on the subsequence, this implies that the whole
sequence $(\nabla u_n)$ converges to $\nabla u$ weakly in
$L^2(\Om;\R^2)$. Taking again $u_n-g_n$ and $u-g$ as test functions
in the equations satisfied
by $u_n$ and $u$, we obtain
$$
\int_{\Om}|\nabla u_n|^2dx=\int_{\Om}\nabla u_n \nabla g_n\,dx\,,
\qquad
\int_{\Om}|\nabla u|^2dx=\int_{\Om}\nabla u \,\nabla g\,dx\,.
$$
As $\nabla u_n \wto \nabla u$ weakly in $L^2(\Om,\R^2)$ and $\nabla
g_n\to\nabla g$ strongly  in $L^2(\Om,\R^2)$, from the
previous equalities we obtain that $\|\nabla u_n\|_{L^2(\Om,\R^2)}$
converges to $\|\nabla u\|_{L^2(\Om,\R^2)}$, which implies the strong
convergence of the gradients in $L^2(\Om,\R^2)$.

By the uniqueness of the gradients of the solutions,
to prove (\ref{2uu*}) it is enough to show that $u^*$ is a solution of
(\ref{**}). This will be
obtained by using Theorem~\ref{due}. First of all we
note that $K\in\Kml$ by Corollary~\ref{Golab2}, and therefore $K$ is
locally connected (see, e.g., \cite[Lemma 1]{Ch-D}).

Let us fix $x\in \overline\Om$ and an open rectangle
$V$ containing $x$. If $x\in\Om$, we  assume that $V\subset\Om$.
If $x\in\partial\Om$, we assume that $V$ is as in the definition of
the Lipschitz part of the boundary. Let $U$ be an open neighbourhood
of $x$ in $\R^2$ such that $U\subset\subset V$. We will prove that
there exists a function
$v\in H^1({U\cap\Om})$, with
$\nabla v=R\,\nabla u^*$ a.e.\ in ${U\cap\Om}$,
such that $v$ is constant q.e.\ on
each connected component of ${U\cap K}$ and of
${U\cap\partial_N\Om}$.
By Theorem~\ref{due} this implies that $u^*$ satisfies~(\ref{**}).

Let $\delta:={\rm dist}(U,\partial V)$. Let us prove that there are
at most $m+{\lambda}/{\delta}$ connected components $C$ of
${\overline V\cap K_n}$ which meet ${U\cap K_n}$.
Indeed, if $C$ meets also
$\partial V$, then $\huno(C)\ge\delta$ (since $C$ connects
a point in $U$ with a point in $\partial V$), so that the number of
these components can not exceed ${\lambda}/{\delta}$. On the
other hand, it is easy to see that the other connected components
of ${\overline V\cap K_n}$ are also connected components of $K_n$,
thus their number can not exceed $m$.

Let $K^1_n,\ldots,K^{k_n}_n$ be the connected components
of ${\overline V\cap K_n}$ which meet ${U\cap K_n}$. As
$k_n\le m+{\lambda}/{\delta}$,
passing to a subsequence we may assume that  $k_n=k$ for every $n$,
and that $K{}^1_n\to \widehat K{}^1$, $\ldots$,
$K{}^k_n\to \widehat K{}^k$
in the Hausdorff metric, where $\widehat K{}^1, \ldots, \widehat K{}^k$
are compact and connected. Arguing as in the proof of Corollary~\ref{Golab2},
we obtain that
\begin{equation}\label{UKsubset}
U\cap K\subset \widehat K{}^1 \cup \cdots\cup \widehat K{}^k\,.
\end{equation}

Let $v_n$ be the harmonic conjugate of $u_n$ in ${V\cap\Om}$ given by
Theorem~\ref{conjug}.
Then $\nabla v_n=R\,\nabla u_n$ a.e.\ in ${V\cap\Om}$.
We may assume that $\int_{V\cap\Om} v_n\,dx=0$ for every $n$.
Since $\nabla v_n=R\,\nabla u_n$ a.e.\ on ${V\cap\Om}$,
we deduce that $(\nabla v_n)$ converges to $R\,\nabla u^*$ weakly in
$L^2({V\cap\Om};\R^2)$, and by the Poincar\'e inequality
$(v_n)$  converges weakly in $H^1({V\cap\Om})$ to a
function $v$ which satisfies
$\nabla v=R\,\nabla u^*$ a.e.\ on ${V\cap\Om}$.

Let us
prove that  for every $i=1,\ldots,k$ there exists a constant
$c^i$ such that $v=c^i$ q.e.\ on $\widehat K{}^i$. This is trivial
when $\widehat K{}^i$ reduces to one point.
If $\widehat K{}^i$ has more than one point,
then $\liminf_{n}{\rm diam}(K^i_n)>0$;
since the sets $K^i_n$ are connected, we obtain also
$\liminf_{n}{\rm cap}(K^i_n)>0$. As $v_n=c^i_n$ q.e.\ on $K^i_n$ for
suitable constants $c^i_n$,
using the Poincar\'e inequality (see, e.g., \cite[Corollary 4.5.3]{Zie})
it follows that $(v_n-c^i_n)$ is bounded in $H^1({V\cap\Om})$, hence
the sequence $(c^i_n)$ is bounded, and therefore, passing to a
subsequence, we may assume that $c^i_n\to c^i$ for a suitable
constant $c^i$. Then $(v_n-c^i_n)$ converges to $v-c^i$ weakly in
$H^1({V\cap\Om})$, and by Lemma \ref{Mosco} we conclude that
$v=c^i$ q.e.\ on $\widehat K{}^i$.

By Proposition \ref{dueconnessi}, if $\widehat K{}^i\cap \widehat
K{}^j\neq \emptyset$, then $v$ is constant q.e.\ on
$\widehat K{}^i\cup \widehat K{}^j$. By (\ref{UKsubset}) this implies that
$v$ is constant q.e.\ on each connected component of ${U\cap K}$.

On the other hand, every $v_n$ is constant q.e.\ on each
connected component of ${V\cap\partial_N\Om}$. Since
$v_n\wto v$ weakly in $H^1({V\cap\Om})$, a sequence of convex 
combinations of the functions $v_n$ converges to $v$ strongly in
$H^1({V\cap\Om})$, and we conclude
that $v$ is constant q.e.\ on each connected component of
${V\cap\partial_N\Om}$, hence on each connected component of
${U\cap\partial_N\Om}$.

Therefore $u^*$ satisfies all hypotheses of Theorem~\ref{due}, which
implies that $u^*$ is a solution of problem~(\ref{**}).
\end{proof}

\end{section}

\begin{section}{COMPACT VALUED INCREASING FUNCTIONS}

In this section we consider increasing functions  $K\colon[0,1]\to\K$, i.e.,
we assume that $K(s)\subset K(t)$ for $s<t$. The following proposition 
extends to compact valued increasing functions a well known 
result about the continuity of real valued monotone functions.

\begin{proposition}\label{diff3}
Let $K\colon [0,1]\to\K$ be an increasing function, and let
$K^-\colon (0,1]\to\K$ and $K^+\colon [0,1)\to\K$ be the functions
defined by
\begin{eqnarray}
&\textstyle K^-(t):={\rm cl}\big(\bigcup_{s<t}K(s)\big)
\qquad\hbox{for }0<t\leq1\,,
\label{k_*t}\\
&\textstyle K^+(t):=\bigcap_{s>t}K(s)\qquad\hbox{for }0
\leq t<1\,,\label{k^*t}
\end{eqnarray}
where ${\rm cl}$ denotes the closure. Then
\begin{equation}\label{7.9}
K^-(t)\subset K(t)\subset K^+(t)
\qquad\hbox{for }0< t<1\,.
\end{equation}
Let $\Theta$ be the set of points $t\in(0,1)$ such that  $K^-(t)=K^+(t)$.
Then ${[0,1]\setmeno\Theta}$ is at most countable, and $K(t_n)\to
K(t)$ in the Hausdorf metric for every $t\in\Theta$ and every sequence
$(t_n)$ in $[0,1]$ converging to~$t$.
\end{proposition}

To prove Proposition~\ref{diff3} we use the following result, which extends
another well known property of real valued monotone functions.

\begin{lemma}\label{diff3.5}
Let $K_1$, $K_2\colon [0,1]\to\K$ be two
increasing functions such that
\begin{equation}\label{K12}
K_1(s) \subset K_2(t)\qquad\hbox{and} \qquad
K_2(s) \subset K_1(t)
\end{equation}
for every $s$, $t\in[0,1]$ with $s<t$.
Let $\Theta$ be the set of points $t\in[0,1]$ such that  $K_1(t)=K_2(t)$.
Then ${[0,1]\setmeno\Theta}$ is at most countable.
\end{lemma}

\begin{proof}
For $i=1,2$, consider the functions
$f_i\colon \overline \Om\times[0,1]\to\R$ defined by
$f_i(x,t):={\rm dist}(x,K_i(t))$, with the convention
that ${\rm dist}(x,\emptyset)={\rm diam}(\Om)$. Then the
functions $f_i(\cdot,t)$ are
Lipschitz continuous with constant $1$ for every $t\in[0,1]$, and
the functions $f_i(x,\cdot)$ are non-increasing for every $x\in\overline\Om$.

Let $D$ be a countable dense subset of $\overline\Om$.  For every
$x\in D$ there exists a countable set $N_x\subset[0,1]$ such that
$f_i(x,\cdot)$ are continuous at every point of $[0,1]\setmeno N_x$. 
By (\ref{K12}) we have $f_1(x,s)\ge f_2(x,t)$ and $f_2(x,s)\ge f_1(x,t)$ for
every $x\in\overline\Om$ and every $s$, $t\in[0,1]$ with $s<t$. This
implies that $f_1(x,t)=f_2(x,t)$ for every $x\in D$ and every
$t\in[0,1]\setmeno N_x$. 
Let $N$ be the countable set defined by
$N:=\bigcup_{x\in D}N_x$, and let $t\in[0,1]\setmeno N$.
Then $f_1(x,t)=f_2(x,t)$ for every $x\in D$, and, 
by continuity, for every 
$x\in\overline\Om$, which yields $K_1(t)=K_2(t)$. 
This proves that ${[0,1]\setmeno N}\subset\Theta$.
\end{proof}

\begin{proof}[Proof of Proposition~\ref{diff3}.]
It is clear that $K^+$ and $K^-$ are increasing and satisfy (\ref{K12}).
Therefore ${[0,1]\setmeno\Theta}$ is at most countable by
Lemma~\ref{diff3.5}.

Let us fix $t\in\Theta$ and a sequence $(t_n)$ in $[0,1]$ converging
to $t$. By the Compactness Theorem~\ref{compactness} we may assume
that $K(t_n)$ converges in the Hausdorff metric to a set $K^*$. For
every $s_1$, $s_2\in[0,1]$, with $s_1<t<s_2$, we have $K(s_1)\subset
K(t_n)\subset K(s_2)$ for $n$ large enough, hence $K(s_1)\subset
K^*\subset K(s_2)$. As $K^*$ is closed this implies $K^-(t)\subset
K^*\subset K^+(t)$, therefore $K^*=K(t)$ by (\ref{7.9}) and by the
definition of~$\Theta$.
\end{proof}

The following result is the analogue of the Helly theorem for
compact valued increasing functions.

\begin{theorem}\label{Helly}
Let $(K_n)$ be a sequence of increasing functions from $[0,1]$ into
$\K$. Then there exist a subsequence, still denoted by $(K_n)$, and an
increasing function $K\colon[0,1]\to\K$, such that $K_n(t)\to K(t)$ in
the Hausdorff metric for every $t\in[0,1]$.
\end{theorem}
\begin{proof}
Let $D$ be a countable dense subset of $(0,1)$. Using a diagonal
argument, we find a subsequence,  still denoted by $(K_n)$, and an
increasing function $K\colon D\to\K$, such that $K_n(t)\to K(t)$ in
the Hausdorff metric for every $t\in D$. Let
$K^-\colon (0,1]\to\K$ and $K^+\colon [0,1)\to\K$ be the increasing
functions defined by
\begin{eqnarray*}
& K^-(t):=\displaystyle{\rm cl}
\Big(\bigcup_{s<t,\; s\in D}K(s)\,\Big)
\qquad\hbox{for }0<t\leq1\,,\\
& K^+(t):=\displaystyle\bigcap_{s>t,\;s\in D }K(s)\qquad\hbox{for }0
\leq t<1\,,
\end{eqnarray*}
where ${\rm cl}$ denotes the closure. Let $\Theta$ be the set of points
$t\in[0,1]$ such that  $K^-(t)=K^+(t)$. As $K^-$ and $K^+$ satisfy
(\ref{K12}), by Lemma~\ref{diff3.5} the set
${[0,1]\setmeno\Theta}$ is at most countable.

Since $K^-(t)\subset K(t)\subset K^+(t)$ for every $t\in D$, we have
$K(t)=K^-(t)=K^+(t)$ for every $t\in {\Theta\cap D}$.
For every $t\in {\Theta\setmeno D}$ we define $K(t):=K^-(t)=K^+(t)$.
To prove that $K_n(t)\to K(t)$ for a given $t\in {\Theta\setmeno D}$,
by the Compactness Theorem~\ref{compactness} we may assume that
$K_n(t)$ converges in the Hausdorff metric to a set $K^*$. For
every $s_1$, $s_2\in D$, with $s_1<t<s_2$, by monotonicity we have
$K(s_1)\subset K^*\subset K(s_2)$. As $K^*$ is closed, this implies
$K^-(t)\subset K^*\subset K^+(t)$, therefore
$K_n(t)\to K(t)$ by the definitions of $\Theta$ and $K(t)$.

Since $[0,1]\setmeno(\Theta\cup D)$ is at most countable, by a
diagonal argument we find a further subsequence, still denoted by
$(K_n)$, and a function
$K\colon{[0,1]\setmeno (\Theta\cup D)}\to \K$, such that $K_n(t)\to
K(t)$  in the Hausdorff metric for every $t\in {[0,1]\setmeno(\Theta\cup D)}$.

Therefore $K_n(t)\to K(t)$ for every $t\in[0,1]$, and this implies that
$K$ is increasing on $[0,1]$.
\end{proof}

For every compact set $K$ in $\R^2$ and
every $g\in\hb(\Omk)$ we define
\begin{equation}\label{e}
\E(g,K):=\min_{v\in {\mathcal V}(g,K)}\Big\{\int_{\Om\setminus K}|\nabla
v|^2dx+\huno(K)\Big\}\,,
\end{equation}
where ${\mathcal V}(g,K)$ is the set introduced  in (\ref{vgk}).

Given a Hilbert space $X$, we recall that $AC([0,1];X)$  is the space of
all absolutely continuous functions defined in $[0,1]$ with
values in $X$. For the main properties of these functions we refer, e.g.,
to \cite[Appendix]{Bre}. Given $g\in AC([0,1];X)$, the time derivative of
$g$, which exists a.e.\ in $[0,1]$, is denoted by $\gdot$. It is
well-known that $\gdot$ is a Bochner integrable function with values
in $X$.

The following result will be crucial in the next section.

\begin{theorem}\label{diff5}
Let $m\ge1$, let $g\in AC([0,1];H^1(\Om))$, and let 
$K\colon[0,1]\to\Kmf$ be an
increasing function.  Suppose that the function  $t\mapsto \E(g(t),K(t))$
is absolutely
continuous on $[0,1]$.
Then the  following conditions are equivalent:
\begin{itemize}
\item[(a)]
\hfil$\displaystyle\frac{d}{ds}\E(g(t),K(s))
\Big|\lower1.5ex\hbox{$\scriptstyle s=t$}=0\quad\hbox{ for a.e.\ }
t\in[0,1]\,,$\hfil
\item[(b)] \hfil$\displaystyle\frac{d}{dt}\E(g(t),K(t))=
2(\nabla u(t)|\nabla \gdot(t))\quad \hbox{ for a.e.\ } t\in[0,1]$,\hfil
\end{itemize}
where $u(t)$ is a solution of the minimum problem (\ref{e}) which defines
$\E(g(t),K(t))$, and $(\cdot|\cdot)$ denotes the scalar product
in $L^2(\Om;\R^2)$.
\end{theorem}

To prove Theorem~\ref{diff5} we need the following lemmas.
\begin{lemma}\label{diff1}
Let $K\in\Kf$ and let
$ F\colon H^1(\Om)\to\R$ be defined by $ F(g)=\E(g,K)$ for every $g\in
H^1(\Om)$. Then $ F$ is of class $C^1$ and for every $g,h\in
H^1(\Om)$ we have
\begin{equation}\label{dphi}
d F(g)\,h=2\int_{\Om\setminus K}\nabla u_g\nabla h\,dx\,,
\end{equation}
where $u_g$ is a solution of the minimum  problem (\ref{e}) which
defines $\E(g,K)$.
\end{lemma}
\begin{proof} Since $u_g$ is a solution of problem (\ref{**}) which
satisfies the boundary condition (\ref{**b}), by linearity
for every $t\in\R$ we have $\nabla u_{g+th}=\nabla u_g+t\nabla u_h$
a.e.\ in $\Om$, hence
\begin{eqnarray*}
  & F(g+th)- F(g)=\displaystyle \int_{\Om\setminus K}|\nabla u_g+t\nabla
u_h|^2\,dx-\int_{\Om\setminus K}|\nabla u_g|^2\,dx=\\
&=2t\displaystyle\int_{\Om\setminus K}\nabla u_g\nabla
u_h\,dx+t^2\int_{\Om\setminus K}|\nabla
u_h|^2\,dx=2t\int_{\Om\setminus K}\nabla u_g\nabla h\,dx+
t^2\int_{\Om\setminus K}|\nabla u_h|^2\,dx\,,
\end{eqnarray*}
where the last equality is deduced from (\ref{**}). Dividing by $t$
and letting $t$ tend to $0$ we obtain (\ref{dphi}). The continuity of
$g\mapsto\nabla u_g$ implies that $ F$ is of class $C^1$.
\end{proof}

Let us consider now the case of time dependent compact sets $K(t)$.

\begin{lemma}\label{diff2}
Let $m\ge 1$ and $\lambda\ge 0$, let
$K\colon [0,1]\to\Kml$ be a function, and
let $ F\colon H^1(\Om){\times}[0,1]\to\R$ be defined by
$ F(g,t)=\E(g,K(t))$. Then the differential $d_1 F$ of $ F$ with
respect to $g$ is continuous at every
point  $(g,t)\in H^1(\Om){\times}[0,1]$ such that
$K(s)\to K(t)$ in the Hausdorff metric as $s\to t$.
\end{lemma}
\begin{proof}
It is enough to apply Lemma~\ref{diff1} and Theorem~\ref{convsol}.
\end{proof}

To deal with the dependence on $t$ of both arguments we need the
following result.

\begin{lemma}\label{diff4}
Let $X$ be a Hilbert space, let $g\in AC([0,1];X)$, and let
$ F\colon X{\times}[0,1]\to\R$ be a function such that $ F(\cdot,t)\in
C^1(X)$ for every $t\in[0,1]$, with differential denoted by
$d_1 F(\cdot,t)$. Let $t_0\in[0,1]$, let $\psi(t):= F(g(t),t)$,
and let $\psi_0(t):= F(g(t_0),t)$. Assume that $t_{0}$ is a
differentiability point of $\psi$ and $g$ and a Lebesgue point of
$\gdot$,
  and that $d_1 F$ is continuous at
$(g(t_0),t_0)$. Then $\psi_0$ is differentiable at $t_0$ and
$$
\dot\psi_0(t_0)=\dot\psi(t_0)-d_1\! F(g(t_0),t_0)\,\gdot(t_0)\,.
$$
\end{lemma}
\begin{proof}
For every $t\in[0,1]$ we have
\begin{eqnarray*}
\psi_0(t)-\psi_0(t_0)=&\hskip-0.5em F(g(t_0),t)- F(g(t),t)+\psi(t)-\psi(t_0)=\\
=&\hskip-1em \displaystyle\int_t^{t_0}d_1 
F(g(s),t)\,\gdot(s)\,ds+\psi(t)-\psi(t_0)\,.
\end{eqnarray*}
The conclusion follows dividing by $t-t_0$ and taking the limit as
$t\to t_0$.
\end{proof}

\begin{proof}[Proof of Theorem~\ref{diff5}.]
Let $ F\colon H^1(\Om){\times}[0,1]\to\R$ be defined by
$ F(g,t)=\E(g,K(t))$. By Proposition~\ref{diff3} and Lemma~\ref{diff2}
$d_1 F$ is continuous in $(g,t)$ for a.e.\ $t\in[0,1]$ and every $g\in
H^1(\Om)$. By Lemmas~\ref{diff1} and~\ref{diff4}
\begin{equation}\label{formula1}
\frac{d}{ds}\E(g(t),K(s))\Big|\lower1.5ex\hbox{$\scriptstyle s=t$}=
\frac{d}{dt}\E(g(t),K(t))-2(\nabla
u(t)|\nabla \gdot(t))\quad\hbox{ for a.e.\ }t\in[0,1]\,.
\end{equation}
The equivalence between (a) and (b) is now obvious.
\end{proof}
\end{section}

\begin{section}{IRREVERSIBLE QUASI-STATIC EVOLUTION}\label{irrev}

In this section we prove the main result of the paper.
\begin{theorem}\label{kt}
Let $m\ge1$, let $g\in AC([0,1];H^1(\Om))$, and let $K_{0}\in\Kmf$. Then
there exists a function $K\colon[0,1]\to\Kmf$ such that
\smallskip
\begin{itemize}
\item[(a)]
\hfil $\displaystyle \vphantom{\frac{d}{ds}}
K_0\subset K(s)\subset K(t)$  for $0\le s\le t\le 1$, \hfil
\item[(b)] \hfil $\displaystyle \vphantom{\frac{d}{ds}}
\E(g(0),K(0))\leq \E(g(0),K)
\quad\forall\, K\in\Kmf,\,\  K\supset K_0$,\hfil
\item[(c)] \hfil $\displaystyle \vphantom{\frac{d}{ds}}
\hbox{for  }\, 0\le t \le1\quad\E(g(t),K(t))\leq \E(g(t),K)
\quad\forall \, K\in\Kmf,\,\  K\supset K(t)$,\hfil
\item[(d)]\hfil $\displaystyle \vphantom{\frac{d}{ds}}
t\mapsto \E(g(t),K(t)) \hbox{ is
absolutely continuous on }[0,1]$, \hfil
\item[(e)]\hfil$\displaystyle\frac{d}{ds}\E(g(t),K(s))
\Big|\lower1.5ex\hbox{$\scriptstyle s=t$}=0\quad
\hbox{for a.e.\ }t\in[0,1]$.\hfil
\end{itemize}
\smallskip
Moreover every function $K\colon[0,1]\to\Kmf$ which satisfies (a)--(e)
satisfies also
\begin{itemize}
\item[(f)] \hfil$\displaystyle\frac{d}{dt}\E(g(t),K(t))=
2(\nabla u(t)|\nabla \gdot(t))\quad \hbox{ for a.e.\ } t\in[0,1]$,\hfil
\end{itemize}
where $u(t)$ is a solution of the minimum  problem (\ref{e}) which defines
$\E(g(t),K(t))$.
\end{theorem}
Here and in the rest of the section $(\cdot|\cdot)$  and $\|\cdot\|$
denote the scalar product  and the norm in $L^2(\Om;\R^2)$.

Theorem~\ref{kt} will be proved by a time discretization process.
Given $\delta>0$, let $N_\delta$ be the largest integer such that
$\delta N_\delta\le 1$;
for $i\geq 0$ let $t_i^\delta:=i\delta$ and, for $0\le i\le
N_\delta$, let $g_i^\delta:=g(t_i^\delta)$.
We define $K_i^\delta$, inductively, as
a solution of the minimum problem
\begin{equation}\label{pidelta}
\min_K\big\{ \E(g_i^\delta,K) : K\in\Kmf,\ K\supset K_{i-1}^\delta\big\}\,,
\end{equation}
where we set $K_{-1}^\delta:=K_0$.



\begin{lemma}
There exists a solution of the minimum problem (\ref{pidelta}).
\end{lemma}
\begin{proof}
By hypothesis $K_{-1}^\delta:=K_0\in\Kmf$. Assume by induction that
$K_{i-1}^\delta\in \Kmf$ and let $\lambda$ be a constant such that
$\lambda>\E(g_i^\delta,K_{i-1}^\delta)$.
Consider a minimizing sequence $(K_n)$ of problem~(\ref{pidelta}).
We may assume that $K_n\in\Kml$ for every $n$.
By the Compactness Theorem~\ref{compactness}, passing to a subsequence,
we may assume that  $(K_n)$ converges in the Hausdorff metric  to
some  compact set $K$ containing $K_{i-1}^\delta$.
For every $n$ let $u_n$ be a
solution of the minimum  problem (\ref{e}) which defines
$\E(g_i^\delta,K_n)$.
By Theorem~\ref{convsol} $(\nabla u_n)$ converges
strongly in $L^2(\Om;\R^2)$ to $\nabla u$, where $u$ is a  solution of
the minimum problem (\ref{e}) which defines $\E(g_i^\delta,K)$.
By Corollary~\ref{Golab2} we have $K\in\Km$ and
$\huno(K)\leq\liminf_{n}\huno(K_n)\le\lambda$, hence $K\in\Kml$. As
$\|\nabla u\|=\lim_{n}\|\nabla u_n\|$,
we conclude that $\E(g_i^\delta,K)\leq\liminf_{n}\E(g_i^\delta,K_n)$.  Since
$(K_n)$ is a minimizing sequence, this proves that $K$  is a
solution of the minimum problem (\ref{pidelta}).
\end{proof}

We define now the step functions $g_\delta$, $K_\delta$, and $u_\delta$
on $[0,1]$
by setting $g_\delta(t):=g_{i}^\delta$, $K_\delta(t):=K_{i}^\delta$,
and $u_\delta(t):=u_{i}^\delta$ for
$t_i^\delta\leq t<t_{i+1}^\delta$, where $u_i^\delta$ is a
solution of the minimum problem (\ref{e}) which defines
$\E(g_{i}^\delta,K_{i}^\delta)$.

\begin{lemma}\label{discr}
There exists a positive function $\rho(\delta)$, converging to zero
as $\delta\to0$, such that
\begin{equation}\label{2discr}
\|\nabla u_j^\delta\|^2+\huno(K_j^\delta)\leq
\|\nabla u_i^\delta\|^2+\huno(K_i^\delta)+
2\int_{t_i^{\delta}}^{t_j^{\delta}}(\nabla u_\delta(t)|\nabla
\gdot(t))\,dt+\rho(\delta)
\end{equation}
for $0\leq i<j\leq N_\delta$.
\end{lemma}
\begin{proof}
Let us fix an integer $r$ with $i\leq r<j$. {}From the absolute
continuity of $g$ we have
$$
g_{r+1}^\delta-g_r^\delta=\int_{t_r^{\delta}}^{t_{r+1}^{\delta}}\gdot(t)\,dt\,,
$$
where the integral is a Bochner integral for functions with values in
$H^1(\Om)$. This implies that
\begin{equation}\label{nabla}
\nabla g_{r+1}^\delta-\nabla g_r^\delta=\int_{t_r^{\delta}}^{t_{r+1}^{\delta}}
\nabla \gdot(t)\,dt\,,
\end{equation}
where the integral is a Bochner integral for functions with values in
$L^2(\Om;\R^2)$.

As $u_r^{\delta}+g_{r+1}^\delta-g_r^\delta\in \hb(\Omk_{r}^\delta)$
and $u_r^{\delta}+g_{r+1}^\delta-g_r^\delta=g_{r+1}^\delta$ q.e.\ on
$\partial_D\Omk_r^\delta$, we have
\begin{equation}\label{a}
\E(g_{r+1}^\delta,K_r^\delta)\le 
\|\nabla u_r^\delta+\nabla g_{r+1}^\delta-\nabla g_r^\delta\|^2+
\huno(K_r^\delta)\,.
\end{equation}
By the minimality of $u_{r+1}^\delta$ and by (\ref{pidelta})
we have
\begin{equation}\label{b}
\|\nabla u_{r+1}^\delta\|^2+\huno(K_{r+1}^\delta)=
\E(g_{r+1}^\delta,K_{r+1}^\delta)\le 
\E(g_{r+1}^\delta,K_r^\delta)\,.
\end{equation}
{}From (\ref{nabla}), (\ref{a}), and (\ref{b})  we obtain
\begin{eqnarray*}
& \|\nabla u_{r+1}^\delta\|^2+\huno(K_{r+1}^\delta)\leq
\|\nabla u_r^\delta+\nabla g_{r+1}^\delta-\nabla g_r^\delta\|^2+
\huno(K_r^\delta)\leq\\
&\leq\|\nabla u_{r}^\delta\|^2+\huno(K_r^\delta)+2
\displaystyle\int_{t_r^{\delta}}^{t_{r+1}^{\delta}}(\nabla u_{r}^\delta|\nabla
\gdot(t))\,dt+\Big(\displaystyle\int_{t_r^{\delta}}^{t_{r+1}^{\delta}}\|\nabla
\gdot(t)\|\,dt\Big)^2\leq\\
&\leq\|\nabla u_{r}^\delta\|^2+\huno(K_r^\delta)+2
\displaystyle\int_{t_r^{\delta}}^{t_{r+1}^{\delta}}
(\nabla u_\delta(t)|\nabla\gdot(t))\,dt+
\sigma(\delta)\displaystyle\int_{t_r^{\delta}}^{t_{r+1}^{\delta}}\|\nabla 
\gdot(t)\|\,dt\,,
\end{eqnarray*}
where
$$
\sigma(\delta):=\max_{i\leq r<j} \int_{t_r^{\delta}}^{t_{r+1}^{\delta}}\|\nabla
\gdot(t)\|\,dt \ \longrightarrow \ 0
$$
by the absolute continuity of the integral.
Iterating now this inequality for $i\leq r<j$ we get
(\ref{2discr}) with $\rho(\delta):=\sigma(\delta)\int_0^1\|\nabla
\gdot(t)\|\,dt$.
\end{proof}

\begin{lemma}\label{estim}
There exists a	constant $\lambda$, depending only on $g$ and $K_0$,
such that
\begin{equation}\label{stima} 
\|\nabla u_i^\delta\|\leq \lambda\quad\hbox{and}\quad
\huno(K_i^{\delta})\leq \lambda
\end{equation}
for every $\delta>0$ and for every  $0\leq i\leq N_\delta$.
\end{lemma}
\begin{proof}
As $g_i^\delta$ is admissible for the problem
(\ref{e}) which defines $\E(g_i^\delta, K_i^\delta)$,
by the minimality of $u_i^\delta$ we have
$\|\nabla u_i^\delta\|\leq\|\nabla g_i^\delta\|$, hence
$\|\nabla u_\delta(t)\|\leq\|\nabla g_\delta(t)\|$ for every
$t\in[0,1]$.
As $t\mapsto g(t)$ is absolutely continuous with values in $H^1(\Om)$
the function $t\mapsto \|\nabla\gdot(t)\|$ is integrable on $[0,1]$
and
there exists a constant $C>0$ such that
$\|\nabla g(t)\|\leq C$ for every $t\in[0,1]$.  This implies the
former
inequality in (\ref{stima}). The latter inequality follows
now from Lemma~\ref{discr} and from the inequality
$\|\nabla u_0^\delta\|^2 +\huno(K_0^\delta)\le 
\|\nabla g(0)\|^2 +\huno(K_0)$, which is an obvious consequence of
(\ref{pidelta}) for ${i=0}$.
\end{proof}

\begin{lemma}\label{Helly2} Let $\lambda$ be the constant of 
Lemma~\ref{estim}.
There exists an increasing function $K\colon[0,1]\to\Kml$ such that,
for every $t\in[0,1]$, $K_{\delta}(t)$ converges to $K(t)$ in the
Hausdorff metric as $\delta\to0$ along a suitable sequence 
independent of~$t$.
\end{lemma}
\begin{proof}
By Theorem \ref{Helly} there exists an increasing function
$K\colon[0,1]\to\K$ such that,
for every $t\in[0,1]$, $K_{\delta}(t)$ converges to $K(t)$ in the
Hausdorff metric as $\delta\to0$ along a suitable sequence
independent of~$t$. By
Lemma~\ref{estim} we have $\huno(K_{\delta}(t))\le\lambda$ for every
$t\in[0,1]$ and every $\delta>0$. By Corollary~\ref{Golab2} this
implies $K(t)\in\Kml$ for every $t\in[0,1]$.
\end{proof}

In the rest of this section, when
we write $\delta\to0$, we always refer to the sequence given
by Lemma~\ref{Helly2}.

For every $t\in[0,1]$ let $u(t)$ be  a solution of  the minimum problem
(\ref{e}) which defines $\E(g(t),K(t))$.
\begin{lemma}\label{l2}
For every $t\in[0,1]$ we have
$\nabla u_{\delta}(t)\to\nabla u(t)$ strongly in $L^2(\Om;\R^2)$.
\end{lemma}
\begin{proof} As $u_{\delta}(t)$ is a solution of the minimum  problem
(\ref{e}) which defines $\E(g_{\delta}(t),K_{\delta}(t))$, and
$g_\delta(t)\to g(t)$ strongly in $H^1(\Om)$, the
conclusion follows from
Theorem~\ref{convsol}.
\end{proof}

\begin{lemma}\label{condb}For every $t\in[0,1]$
we have
\begin{equation}\label{pt}
\E(g(t),K(t))\leq \E(g(t),K)\quad\forall\,K\in\Kmf\,,\ K\supset
K(t)\,.
\end{equation}
Moreover
\begin{equation}\label{pt0}
\E(g(0),K(0))\leq \E(g(0),K)\quad\forall\,K\in\Kmf\,,\ K\supset
K_0\,.
\end{equation}
\end{lemma}
\begin{proof} Let us fix $t\in[0,1]$ and $K\in\Kmf$ with $K\supset K(t)$.
Since  $K_\delta(t)$ converges
to $K(t)$ in the Hausdorff metric as $\delta\to0$,
by Lemma~\ref{differ} there exists a sequence
$(K_\delta)$ in $\Kmf$,
converging to $K$ in the Hausdorff metric, such that
$K_\delta\supset K_\delta(t)$ and
$\huno(K_\delta\setmeno K_\delta(t))\to \huno(K\setmeno K(t))$
as $\delta\to0$.
By Lemma \ref{estim} this implies that $\huno(K_\delta)$ is bounded as 
$\delta\to0$.

Let $v_\delta$ and $v$ be solutions of the minimum problems
(\ref{e}) which define $\E(g_\delta(t),K_\delta)$ and
$\E(g(t),K)$, respectively. By Theorem~\ref{convsol}
$\nabla v_\delta\to\nabla v$ strongly
in $L^2(\Om;\R^2)$.
The minimality of $K_\delta(t)$ expressed by (\ref{pidelta})
gives
$\E(g_\delta(t),K_\delta(t))\leq \E(g_\delta(t),K_\delta)$,
which implies
$\|\nabla u_\delta(t)\|^2\leq\|\nabla v_\delta\|^2+\huno(K_\delta\setmeno
K_\delta(t))$.
Passing to the limit as $\delta\to0$ and using Lemma~\ref{l2}
we get $\|\nabla u(t)\|^2\leq\|\nabla v\|^2+\huno(K\setmeno K(t))$.
Adding $\huno(K(t))$ to both sides we obtain (\ref{pt}).

A similar proof holds for (\ref{pt0}).
By (\ref{pidelta}) we have
$\E(g_\delta(0),K_\delta(0)) \le \E(g_\delta(0),K)=\E(g(0),K)$, which implies
$\|\nabla u_\delta(0)\|^2+\huno(K_\delta(0))\leq \E(g(0),K)$.
Passing to the limit as $\delta\to0$ and using Lemma~\ref{l2}
and Corollary \ref{Golab2} we obtain~(\ref{pt0}).
\end{proof}

The previous lemma proves conditions (b) and (c) of Theorem~\ref{kt}.
To show that  conditions (d) and (e) are also satisfied,
we begin by proving the following inequality.
\begin{lemma}\label{ineq}
For every $s,\, t$ with $0\leq s<t\leq1$
\begin{equation}\label{diseg}
\|\nabla u(t)\|^2+\huno(K(t))\leq\|\nabla u(s)\|^2+\huno(K(s))+
2\int_s^t(\nabla u(\tau)|\nabla \gdot(\tau))d\tau\,.
\end{equation}
\end{lemma}
\begin{proof}
Let us fix $s,t$ with $0\leq s<t\leq1$. Given $\delta>0$ let $i$ and $j$
be the
integers such that $t_i^\delta\leq s<t_{i+1}^\delta$ and
$t_j^\delta\leq t<t_{j+1}^\delta$.
Let us define $s_\delta:=t_i^\delta$ and $t_\delta:=t_j^\delta$.
Applying Lemma~\ref{discr} we obtain
\begin{equation}\label{ediscr}
\|\nabla u_\delta(t)\|^2+\huno(K_\delta(t)\setmeno K_\delta(s))\leq
\|\nabla u_\delta(s)\|^2+
2\int_{s_{\delta}}^{t_{\delta}}\!\!\!(\nabla u_\delta(\tau)|\nabla
\gdot(\tau))\,d\tau+\rho(\delta)\,,
\end{equation}
with $\rho(\delta)$ converging to zero as $\delta\to0$. 
By Lemma~\ref{l2}
for every $\tau\in[0,1]$ we have $\nabla u_{\delta}(\tau)\to \nabla u(\tau)$
strongly in
$L^2(\Om,\R^2)$ as $\delta\to0$, and by Lemma~\ref{estim} we have
$\|\nabla u_\delta(\tau)\|\le \lambda$ for every $\tau\in[0,1]$.
By Corollary~\ref{sci2} we get
$$
\huno(K(t)\setmeno K(s))\leq\liminf_{\delta\to0}\, 
\huno(K_\delta(t)\setmeno K_\delta(s))\,.
$$
Passing now to the limit in (\ref{ediscr}) as $\delta\to0$ we
obtain (\ref{diseg}).
\end{proof}

The following lemma concludes the proof of Theorem~\ref{kt},
showing that also conditions (d), (e) and (f) are satisfied.

\begin{lemma}
The function $t\mapsto \E(g(t),K(t))$ is absolutely continuous on
$[0,1]$ and
\begin{equation}\label{gtkt}
\frac{d}{dt}\E(g(t),K(t))=2(\nabla u(t)|\nabla\gdot(t))\qquad\hbox{ for
a.e.\ }t\in[0,1]\,.
\end{equation}
Moreover
\begin{equation}\label{gtks}
\frac{d}{ds}\E(g(t),K(s))\Big|\lower1.5ex\hbox{$\scriptstyle 
s=t$}=0\qquad\hbox{ for a.e.\ }t\in[0,1]\,.
\end{equation}
\end{lemma}
\begin{proof}Let $0\leq s<t\leq1$.
{}From the previous lemma we get
\begin{equation}\label{alto}
\E(g(t),K(t))-\E(g(s),K(s))\leq2\int_s^t(\nabla
u(\tau)|\nabla\gdot(\tau))\,d\tau\,.
\end{equation}
On the other hand, by condition (c) of Theorem~\ref{kt} we have
$\E(g(s),K(s))\leq \E(g(s),K(t))$, and by Lemma~\ref{diff1}
$$
\E(g(t),K(t))-\E(g(s),K(t))=
2\int_s^t(\nabla u(\tau, t)|\nabla \gdot(\tau))\,d\tau\,,
$$
where $u(\tau,t)$ is a solution of the minimum problem~(\ref{e}) which defines
$\E(g(\tau),K(t))$.
Therefore
\begin{equation}\label{basso}
\E(g(t),K(t))-\E(g(s),K(s))\geq
2\int_s^t(\nabla u(\tau, t)|\nabla \gdot(\tau))\,d\tau\,.
\end{equation}
Since there exists a constant $C$ such that
$\|\nabla u(\tau)\|\leq\|\nabla g(\tau)\|\leq C$ and
$\|\nabla u(\tau,t)\|\leq\|\nabla g(\tau)\|\leq C$ for $s\leq\tau\leq t$,
from (\ref{alto}) and (\ref{basso})
we obtain
$$
\big|\E(g(t),K(t))-\E(g(s),K(s))\big|\leq 2\,C
\int_s^t\|\nabla \gdot(\tau)\|\,d\tau\,,
$$
which proves that the  function $t\mapsto \E(g(t),K(t))$ is absolutely
continuous.

As $\nabla u(\tau,t)\to\nabla u(t)$ strongly in $L^2(\Om,\R^2)$ when
$\tau\to t$, if we divide (\ref{alto}) and (\ref{basso}) by $t-s$, and
take the limit as $s\to t-$ we obtain (\ref{gtkt}).
Equality (\ref{gtks}) follows from Theorem~\ref{diff5}.
\end{proof}

Theorem~\ref{kt0} is a consequence of Theorem~\ref{kt} and of
the following lemma.
\begin{lemma}\label{previous}
Let $K\colon[0,1]\to\Kmf$ be a function which satisfies
conditions (a)--(e) of Theorem~\ref{kt}. Then, for $0<t\le1$,
\begin{equation}\label{s<t}
\E(g(t),K(t))\leq \E(g(t),K)
\quad\forall \,K\in\Kmf\,\ K\supset \textstyle
\bigcup_{s<t}K(s)\,.
\end{equation}
\end{lemma}

\begin{proof}
Let us fix $t$, with $0<t\le1$, and $K\in\Kmf$, with $K\supset
\bigcup_{s<t}K(s)$. For  $0\leq s<t$ we have $K\supset K(s)$, and from
condition (c) of Theorem~\ref{kt} we obtain $\E(g(s),K(s))\leq
\E(g(s),K)$. As the functions $s\mapsto \E(g(s),K(s))$ and
$s\mapsto \E(g(s),K)$ are continuous, passing to the limit as $s\to
t-$ we get (\ref{s<t}).
\end{proof}

The following lemma shows that $K(t)$, $K^-(t)$, and $K^+(t)$ have the
same total energy.
\begin{lemma}\label{*****}
Let $K\colon[0,1]\to\Kmf$ be a function which satisfies
conditions (a)--(e) of Theorem~\ref{kt}, and let $K^-(t)$ and
$K^+(t)$ be defined by (\ref{k_*t}) and (\ref{k^*t}). Then
\begin{eqnarray}
&&\E(g(t),K(t))=\E(g(t),K^-(t))\quad\hbox{for }0<t\le1\,,\label{eq*}\\
&&\E(g(t),K(t))= \E(g(t),K^+(t))\quad\hbox{for }0\le t<1\,.\label{eq**}
\end{eqnarray}
\end{lemma}

\begin{proof} Let $0<t\le1$. Since $K(s)\to K^-(t)$ in the Hausdorff
metric as $s\to t-$, and $\huno(K(s))\to \huno(K^-(t))$ by
Corollary~\ref{Golab2}, it follows that $\E(g(s),K(s))\to \E(g(t),K^-(t))$ as
$s\to
t-$  by Theorem~\ref{convsol}. As the function $s\mapsto \E(g(s),K(s))$
is continuous, we obtain (\ref{eq*}). The proof of (\ref{eq**}) is
analogous.
\end{proof}

\begin{remark}\label{rem}
{\rm {}From Lemmas~\ref{previous} and~\ref{*****} it
follows that, if $K\colon[0,1]\to\Kmf$ is a function which satisfies
conditions (a)--(e) of Theorem~\ref{kt},  the same is true for the
functions
$$
t\mapsto\left\{\begin{array}{ll}
\hskip-0.5em K(0) & \hbox{for }t=0\,,\\
\hskip-0.5em K^-(t)& \hbox{for } 0<t\leq 1
\,,
\end{array}\right.\qquad
t\mapsto\left\{\begin{array}{ll}
\hskip-0.5em K^+(t)& \hbox{for }0\leq t<1\,,\\
\hskip-0.5em K(1)& \hbox{for } t=1
\,,
\end{array}\right.
$$
where $K^-(t)$ and
$K^+(t)$ are defined by (\ref{k_*t}) and (\ref{k^*t}). Therefore the
problem has a left-continuous  solution and a
right-continuous solution.
}\end{remark}

\begin{remark}\label{g0}
{\rm In Theorem~\ref{kt} suppose that
$\E(g(0),K_0)\le \E(g(0),K)$
for every $K\in\Kmf$ with $K\supset K_0$.
Then in our time discretization process we can take
$K_0^\delta=K_0$ for every $\delta>0$.
Therefore there exists a function $K\colon[0,1]\to\Kmf$,
satisfying conditions (a)--(e) of Theorem~\ref{kt}, such that
$K(0)=K_0$. In particular this happens for every  $K_0$
whenever $g(0)=0$.

In the case $g(0)=0$,  by condition (b) we must
have $\huno({K(0)\setmeno K_0})=0$, hence
$K(0)=K_0$ if $K(0)$ has no isolated points.
If we disregard this natural constraint and $K_0$ has $m_0$
connected components, for every finite set
$F\subset\overline\Om$ with no more than $m-m_0$ elements
we can find also a solution with
$K(0)={K_0\cup F}$. Indeed in our time discretization process we can take
$K_0^\delta={K_0\cup F}$ for every $\delta>0$.
}\end{remark}

We consider now the case where $g(t)$ is proportional to a fixed
function $h\in H^1(\Om)$.

\begin{proposition}\label{th}
In Theorem~\ref{kt} suppose that $g(t)=\varphi(t)\,h$, where
$\varphi\in AC([0,1])$ is non-decreasing and non-negative, and $h$ is
a fixed function in $H^1(\Om)$. Let $K\colon[0,1]\to\Kmf$ be a 
function which satisfies
conditions (a)--(e) of Theorem~\ref{kt}. Then
\begin{equation}
\E(g(t),K(t))\le \E(g(t),K(s))
\end{equation}
for $0\le s<t\le1$.
\end{proposition}

\begin{proof} Let us fix $0\le s<t\le1$. For every $\tau\in[0,1]$ let
$v(\tau)$ be a solution of the minimum problem (\ref{e}) which defines
$\E(h,K(\tau))$. As $u(\tau)=\varphi(\tau)\,v(\tau)$ and
$\gdot(\tau)=\dot\varphi(\tau)\,h$,
from condition (f) we obtain, adding and subtracting
$\E(g(s),K(s))$,
\begin{eqnarray*}
&\E(g(t),K(t))-\E(g(t),K(s))=\\
&\displaystyle=2\int_s^t(\nabla v(\tau)|\nabla h)
\,\varphi(\tau)\,\dot{\varphi}(\tau)\,d\tau +
(\varphi(s)^2-\varphi(t)^2)\|\nabla v(s)\|^2\,.
\end{eqnarray*}
As $v(\tau)$ is a solution of problem (\ref{**})
with $K=K(\tau)$, and
$v(\tau)=h$ q.e.\ on $\partial_D\Omk(\tau)$, we have
$(\nabla v(\tau)|\nabla h)=\|\nabla v(\tau)\|^2$.
By the monotonicity of $\tau\mapsto K(\tau)$, for $s\le \tau\le t$
we have $v(s)\in \hb(\Omk(\tau))$ and $v(s)=h$
q.e.\ on $\partial_D\Omk(\tau)$. By the minimum property of $v(\tau)$
we obtain $\|\nabla v(\tau)\|^2\le \|\nabla v(s)\|^2$ for $s\le
\tau\le t$.
Therefore
\begin{eqnarray*}
&\E(g(t),K(t))-\E(g(t),K(s))\le\\
&\displaystyle
\le 2 \int_s^t \varphi(\tau)\,\dot{\varphi}(\tau)\, d\tau \;
\|\nabla v(s)\|^2 +
(\varphi(s)^2-\varphi(t)^2)\|\nabla v(s)\|^2= 0\,,
\end{eqnarray*}
which concludes the proof.
\end{proof}
\end{section}

\begin{section}{BEHAVIOUR NEAR THE TIPS}\label{tips}

In this section we consider a function $K\colon[0,1]\to\Kmf$
which satisfies conditions (a)--(e) of Theorem~\ref{kt} for a suitable $g\in
AC([0,1];H^1(\Om))$, and we study the behaviour of
the solutions $u(t)$ near the ``tips'' of the sets $K(t)$. Under
some natural assumptions, we shall see that $K(t)$ satisfies Griffith's
criterion for crack growth.

For every bounded open set $A\subset\R^2$ with Lipschitz boundary,
for every compact set $K\subset\R^2$, and for every function
$g\colon \partial A\setmeno K\to\R$ we define
\begin{equation}\label{eloc}
\E(g,K, A):=\min_{v\in{\mathcal V}(g,K,A)}
\Big\{\int_{A\setminus K}|\nabla v|^2\,dx
+\huno(K\cap \overline A)\Big\} \,,
\end{equation}
where
$$
{\mathcal V}(g,K,A):=\{v\in\hb(A\setmeno K): v=g \quad\hbox{q.e.\ on }
\partial A\setmeno K\}\,.
$$

We now consider in particular the case where $K$ is a regular arc, and
summarize some known results on the behaviour of a solution of
problem (\ref{**}) near the end-points of $K$. Let $B$ be an open ball
in $\R^2$ and let $\gamma\colon [\sigma_0, \sigma_1]\to \R^2$ be a
simple path of class $C^2$ parametrized by arc length. Assume that
$\gamma(\sigma_0)\in\partial B$ and $\gamma(\sigma_1)\in\partial B$,
while $\gamma(\sigma)\in B$ for $\sigma_0< \sigma< \sigma_1$. Assume
in addition that $\gamma$ is not tangent to $\partial B$ at $\sigma_0$
and $\sigma_1$.
For every $\sigma\in [\sigma_0, \sigma_1]$ let ${\Gamma}(\sigma):=
\{\gamma(s): \sigma_0\le s \le \sigma\}$.

\begin{theorem}\label{Grisvard1}
Let $\sigma_0< \sigma< \sigma_1$ and let $u$ be a solution to
problem (\ref{**}) with $\Omega=B$,  $\partial_D\Omega=\partial B$,
and $K={\Gamma}(\sigma)$.
Then there exists a unique constant $\kappa=\kappa(u,\sigma)\in\R$ such that
\begin{equation}\label{sif}
u-\kappa \,\sqrt{2\rho/\pi}\,
\sin(\theta/2)\in H^2(B\setmeno {\Gamma}(\sigma)) \cap
H^{1,\infty}(B\setmeno {\Gamma}(\sigma)) \,,
\end{equation}
where $\rho(x)=|x-\gamma(\sigma)|$ and $\theta(x)$ is the continuous 
function on
$B\setmeno {\Gamma}(\sigma)$ which coincides with the oriented angle
between $
\dot\gamma(\sigma)$ and $x-\gamma(\sigma)$, and
vanishes on the points of the form
$x=\gamma(\sigma)+\e\, 
\dot \gamma(\sigma)$ for
sufficiently small $\e>0$.
\end{theorem}

\begin{proof}
Let $B^-$ and $B^+$ be the connected components of
$B\setmeno {\Gamma}(\sigma_1)$. Since $B^-$ and $B^+$ have a Lipschitz
boundary, by Proposition~\ref{h1} $u$ belongs to $H^1(B^-)$ and
$H^1(B^+)$. This implies that $u\in L^2(B)$, and hence
$u\in H^1(B\setmeno {\Gamma}(\sigma))$. The conclusion follows now from
\cite[Theorem 4.4.3.7 and Section 5.2]{Gri1}, as shown in
\cite[Appendix~1]{MSh}.
\end{proof}

\begin{remark}\label{stress}{\rm
If $u$ is interpreted as the third component of the displacement in an
anti-plane shear, as we did in the introduction, then $\kappa$
coincides with the {\it Mode III stress intensity factor\/} $K_{I\!I\!I}$
of the displacement $(0,0,u)$.
}\end{remark}

\begin{theorem}\label{Grisvard2}
Let $g\colon \partial B\setmeno \{\gamma(\sigma_0)\} \to\R$ be a 
function such that for
every $\sigma_0< \sigma< \sigma_1$ there exists
$g(\sigma)\in \hb(B\setmeno {\Gamma}(\sigma))$ with $g(\sigma)=g$
q.e.\ on $\partial B\setmeno {\Gamma}(\sigma)=
\partial B\setmeno \{\gamma(\sigma_0)\}$. Let $v(\sigma)$ be a
solution of the minimum problem (\ref{eloc}) which defines
$\E(g, {\Gamma}(\sigma), B)$.
Then, for every $\sigma_0< \sigma< \sigma_1$,
$$
\frac{d}{d\sigma}\E(g, {\Gamma}(\sigma), B)=
1-\kappa(v(\sigma),\sigma)^2\,,
$$
where $\kappa$ is defined by (\ref{sif}).
\end{theorem}

\begin{proof}
It is enough to adapt the proof of \cite[Theorem 6.4.1]{Gri2}.
\end{proof}

Let us return to the function $K\colon[0,1]\to\Kmf$
considered at the beginning of the section, and let
$0\le t_0<t_1\le 1$. Suppose that the following structure condition is
satisfied: there exist a finite family of simple arcs ${\Gamma}_i$,
$i=1,\ldots,p$, contained in $\Om$ and parametrized by arc length
by $C^2$  paths $\gamma_i\colon [\sigma_i^0, \sigma_i^1]\to \Om$,
such that, for $t_0<t<t_1$,
\begin{equation}\label{structure}
K(t)=K(t_0)\cup \bigcup_{i=1}^p {\Gamma}_i(\sigma_i(t))\,,
\end{equation}
where ${\Gamma}_i(\sigma):=
\{\gamma_i(\tau): \sigma_i^0\le \tau \le \sigma\}$ and
$\sigma_i\colon [t_0,t_1]\to [\sigma_i^0, \sigma_i^1]$ are
non-decreasing functions with $\sigma_i(t_0)=\sigma_i^0$ and
$\sigma_i^0<\sigma_i(t)<\sigma_i^1$ for $t_0<t<t_1$.
Assume also that the arcs
${\Gamma}_i$ are pairwise disjoint, and that
${\Gamma}_i\cap K(t_0)=\{\gamma_i(\sigma_i^0)\}$.
We consider the sets ${\Gamma}_i(\sigma_i(t))$ as the increasing
branches of the fracture $K(t)$ and the
points $\gamma_i(\sigma_i(t))$ as their moving tips.
For $i=1,\ldots,p$ and
$\sigma_i^0<\sigma<\sigma_i^1$
let $\kappa_i(u,\sigma)$ be the stress intensity factor
defined by (\ref{sif}) with $\gamma=\gamma_i$ and $B$ equal to a
sufficiently small
ball centred at $\gamma_i(\sigma)$.

We are now in a position to state the main result of this section.

\begin{theorem}\label{Griffith}
Let $m\ge1$, let $K\colon[0,1]\to\Kmf$ be a function which satisfies
conditions (a)--(e) of Theorem~\ref{kt} for a suitable
$g\in AC([0,1];H^1(\Om))$, let $u(t)$ be a solution of the minimum
problem (\ref{e}) which defines $\E(g(t),K(t))$,
and let $0\le t_0<t_1\le 1$.
Assume that (\ref{structure}) is satisfied for $t_0<t<t_1$, and that
the arcs ${\Gamma}_i$ and the functions $\sigma_i$
satisfy all properties considered above.
Then
\begin{eqnarray}
&\dot\sigma_i(t)\ge 0\quad \hbox{ for a.e.\ } t\in (t_0,t_1)\,,
\label{sigmadot}\\
&1-\kappa_i(u(t),\sigma_i(t))^2\ge 0
\quad \hbox{ for every\ } t\in (t_0,t_1)\,,
\label{sif>}\\
&\big\{1- \kappa_i(u(t),\sigma_i(t))^2\big\}
\,\dot\sigma_i(t)= 0\quad \hbox{ for a.e.\ } t\in (t_0,t_1)\,,
\label{sif=}
\end{eqnarray}
for $i=1,\ldots,p$.
\end{theorem}

The first condition says simply that the length of every branch
of the fracture can
not decrease, and reflects the irreversibility of the process.
The second condition says
that the absolute value of the stress intensity factor must be less 
than or equal to
$1$ at each tip and for every time. The last condition says that,
at a given tip, the stress intensity factor must be equal to $\pm 1$
at almost every time in which this tip moves with a positive velocity.
This is Griffith's criterion for crack growth in our model.

To prove Theorem~\ref{Griffith} we use the following lemma.

\begin{lemma}\label{cond3loc} Let $m\ge1$, let $H\in\Kmf$ with
$h$ connected components,
let $g\in H^1(\Om)$, and let $u$ be a solution of the
minimum problem (\ref{e}) which
defines $\E(g,H)$. Given  an open subset $A$ of $\Om$, with
Lipschitz boundary, such that $H\cap\overline A\neq\emptyset$, let
$q$ be the number of connected components of $H$ which meet
$\overline A$.
Assume that
\begin{equation}\label{minn}
\E(g,H)\leq \E(g,K)\qquad\forall \,K\in\Kmf,\,\ K\supset
H\,.
\end{equation}
Then
\begin{equation}\label{minnA}
\E(u,H,A)\leq \E(u,K,A)\qquad\forall \,K\in\KA,\,\ K\supset
H\cap\overline A\,.
\end{equation}
\end{lemma}

\begin{proof}
Let $K\in\KA$ with $K\supset
H\cap\overline A$,  let $v$ be  a solution of the
minimum problem (\ref{eloc})  which defines $\E(u,K,A)$,
and let $w$ be the function defined
by $w:=v$ on $\overline A\setmeno K$ and by $w:=u$ on
$(\overline\Om\setmeno \overline A)\setmeno  H$. As $v=u$ q.e.\ on
$\partial A\setmeno K$ the function $w$ belongs to $\hb(\Om\setmeno
(H\cup K))$; using also the fact that $u=g$ q.e.\ on  $\partial_D\Om\setmeno
H$, we obtain that  $w=g$ q.e.\ on $\partial_D\Om\setmeno (H\cup K)$.
Therefore
\begin{eqnarray}
&\displaystyle
\E(g,H\cup K)\leq\int_{\Om\setminus (H\cup K)}|\nabla
w|^2\,dx+\huno(H\cup K)=\label{pprima}\\
&\displaystyle
= \int_{A\setminus K}|\nabla v|^2\,dx+\huno(K\cap \overline A)+
\int_{(\Om\setminus A)\setminus H}|\nabla u|^2\,dx+
\huno(H\setmeno\overline A)\,.\nonumber
\end{eqnarray}
On the other hand, by the minimality of $u$,
\begin{eqnarray}
&\displaystyle
\int_{A\setminus H}|\nabla u|^2\,dx+\huno(H\cap \overline A)+
\int_{(\Om\setminus A)\setminus H}|\nabla u|^2\,dx+\huno(H\setmeno\overline
A)=\label{sseconda}\\
&\displaystyle
=\int_{\Om\setminus H}|\nabla u|^2\,dx+\huno(H)=\E(g,H)\leq\E(g,H\cup K)\,,
\nonumber
\end{eqnarray}
where the last inequality follows from (\ref{minn}), since $H\cup K$
has no more than $m$  connected components (indeed, ${H\cup K}$ has 
exactly ${h-q}$ connected components which do not meet $\overline A$, 
and every connected component of ${H\cup K}$ which meets $\overline A$
contains a connected component of $K$, so that their number does not 
exceed ${q+m-h}$).
{}From (\ref{pprima}) and (\ref{sseconda}) we obtain
$$
\int_{A\setminus H}|\nabla u|^2\,dx+ \huno(H\cap \overline A)\leq
\int_{A\setminus K}|\nabla v|^2\,dx+ \huno(K\cap \overline A)\,,
$$
and the minimality of $v$ yields (\ref{minnA}).
\end{proof}

\begin{proof}[Proof of Theorem \ref{Griffith}.]
Let $t$ be an arbitrary point in $(t_0,t_1)$ and let
$B_i$,  $i=1,\ldots,p$, be a family of open
balls centred at the points $\gamma_i(\sigma_i(t))$.
If the radii are sufficiently small, we have $\overline
B_i\subset \Om$ and $\overline B_i\cap K(t_0)=
\overline B_i\cap \overline B_j=\overline B_i\cap
{\Gamma}_j=\emptyset$ for $j\neq i$. Moreover we may assume
that $B_i\cap {\Gamma}_i=
\{\gamma_i(\sigma):\tau_i^0<\sigma<\tau_i^1\}$,
for suitable constants $\tau_i^0,\,\tau_i^1$ with
$\sigma_i^0<\tau_i^0<\sigma_i(t)<\tau_i^1<\sigma_i^1$, and that
the arcs ${\Gamma}_i$
intersect $\partial B_i$ only at the points $\gamma_i(\tau_i^0)$ and
$\gamma_i(\tau_i^1)$, with a transversal intersection.
All these properties, together with (\ref{structure}), imply that
\begin{equation}\label{Ga}
\overline B_i\cap K(s)=\overline B_i\cap {\Gamma}_i(\sigma_i(s))=
  \{\gamma_i(\sigma):\tau_i^0\leq\sigma\le\sigma_i(s)\}
\quad\hbox{if}\quad \tau_i^0<\sigma_i(s)<\tau_i^1\,.
\end{equation}
In particular this happens for $s=t$, and for $s$ close to $t$ if $\sigma_i$
is continuous at~$t$.

By condition~(c) of Theorem~\ref{kt}
and by Lemma~\ref{cond3loc} for every $i$ we
have that
$$
\E(u(t),K(t),B_i)\le \E(u(t),K,B_i)
\quad\forall \, K\in{\mathcal K}_1^f(\overline B_i),\,\
K\supset K(t)\cap\overline B_i\,.
$$
By (\ref{Ga}) this implies, taking $K:={\Gamma}_i(\sigma)\cap \overline B_i=
\{\gamma_i(\tau):\tau_i^0\leq\tau\le\sigma\}$,
$$
\E(u(t),{\Gamma}_i(\sigma_i(t)),B_i)\le
\E(u(t),{\Gamma}_i(\sigma),B_i)\quad \hbox{ for } \sigma_i(t)\le
\sigma\le \tau_i^1\,,
$$
which yields
\begin{equation}\label{ge0}
\frac{d}{d\sigma} \E(u(t),{\Gamma}_i(\sigma),B_i)
\Big|\lower1.5ex\hbox{$\scriptstyle \sigma=\sigma_i(t)$}
\ge 0\,.
\end{equation}
Inequality (\ref{sif>}) follows now from Theorem~\ref{Grisvard2} applied
with $g:=u(t)$.

By condition (e) of Theorem \ref{kt} for
a.e.\ in $t\in (t_0,t_1)$ we have
$\frac{d}{ds}\E(g(t),K(s))|_{s=t}=0$. Moreover, for a.e.\ in $t\in (t_0,t_1)$
the derivative $\dot\sigma_i(t)$ exists
for $i=1,\ldots,p$.
Let us fix $t\in (t_0,t_1)$ which satisfies all these properties.

By (\ref{Ga}) for $s$ close to $t$ we have
\begin{equation}\label{subadd}
\E(g(t),K(s))\le \sum_{i=1}^{p}
\E(u(t),{\Gamma}_i(\sigma_i(s)),B_i) + \E(u(t),K,A)\,,
\end{equation}
where $K:=K(t_0)\cup \bigcup_i{\Gamma}_i(\tau_i^0)$ and
$A:=\Om\setmeno \bigcup_i \overline B_i$. Note that the equality holds
in (\ref{subadd}) for $s=t$. As the functions $s\mapsto \E(g(t),K(s))$
and $s\mapsto \E(u(t),{\Gamma}_i(\sigma_i(s)),B_i)$ are differentiable
at $s=t$ (by Theorem~\ref{Grisvard2} and by the existence of
$\dot\sigma_i(t)$), we
conclude that
\begin{eqnarray*}
&\displaystyle
0=\frac{d}{ds} \E(g(t),K(s))
\Big|\lower1.5ex\hbox{$\scriptstyle s=t$}
=\sum_{i=1}^{p}
\frac{d}{ds} \E(u(t),{\Gamma}_i(\sigma_i(s)),B_i)
\Big|\lower1.5ex\hbox{$\scriptstyle s=t$}=\\
&\displaystyle
=\sum_{i=1}^{p} \frac{d}{d\sigma} \E(u(t),{\Gamma}_i(\sigma),B_i)
\Big|\lower1.5ex\hbox{$\scriptstyle \sigma=\sigma_i(t)$}
\, \dot\sigma_i(t)=
\sum_{i=1}^{p} \big\{1- \kappa_i(u(t),\sigma_i(t))^2\big\}
\,\dot\sigma_i(t)\,.
\end{eqnarray*}
By (\ref{sigmadot}) and (\ref{sif>}) we have
  $\big\{1- \kappa_i(u(t),\sigma_i(t))^2\big\}
\,\dot\sigma_i(t)\geq0$ for $i=1,\ldots, p$, so that the previous
equalities yield (\ref{sif=}).
\end{proof}

\end{section}


\end{document}